\documentclass[12pt]{article}

\usepackage{graphicx, epsfig}
\usepackage{amsmath,amssymb}
\usepackage[authoryear,nonamebreak,elide]{natbib}
\usepackage[latin1]{inputenc}
\usepackage{multirow}
\usepackage{caption}
\usepackage[toc,page]{appendix}
\usepackage{epstopdf}
\usepackage{cancel}
\usepackage[normalem]{ulem}

\usepackage{color}
\usepackage{bm}

\newcommand{\dd}{\mathrm{d}}

\usepackage{color}
\usepackage{bibentry} 

\usepackage[para]{threeparttable} 
\usepackage{lscape} 

\usepackage{amsthm} 

\usepackage{verbatim}

\usepackage{alltt}

\allowdisplaybreaks

\begin{document}

\title{Higher-order approximate confidence intervals}
\author{Eliane C. Pinheiro
\ \  \ \ Silvia L.P. Ferrari\thanks{Corresponding author. Email: silviaferrari@usp.br} \\
{\small {\em Department of Statistics, University of S\~ao Paulo,
Brazil}}\\
Francisco M.C. Medeiros \\
{\small {\em Department of Statistics, Federal University of Rio Grande do Norte, Brazil}}
}
\maketitle
\begin{abstract}

Standard confidence intervals employed in applied statistical analysis are usually based on asymptotic approximations. Such approximations can be considerably inaccurate in small and moderate sized samples. We derive accurate confidence intervals based on higher-order approximate quantiles of the score function.
The coverage approximation error is $O(n^{-3/2})$ while the approximation error of confidence intervals based on the  first-order asymptotic distribution of the Wald, score, and signed likelihood ratio statistic is $O(n^{-1/2})$. 
Monte Carlo simulations confirm the theoretical findings. An implementation for regression models and real data applications are provided.


\end{abstract}

\section{Introduction}
\label{introduction}

In applied statistics it is often the case that the practitioner wishes to estimate the parameters of a statistical model fitted to a data set. A point estimate is a single realization from the distribution of the possible values of the chosen estimator and does not carry information on its uncertainty. An uncertainty description is provided by confidence intervals. 

Confidence intervals can, in principle, be constructed from the distribution of the estimator, but the exact distribution is unknown except in particular situations. The most common approach for estimating parameters of a parametric family of distributions is the use of Wald-type confidence intervals, i.e. maximum likelihood estimates (MLEs) coupled with the corresponding first-order asymptotic normal approximation. For large sample sizes, this approach is reasonably reliable and leads to virtually unbiased estimates and confidence intervals with approximately correct coverage. In small and moderate sized samples, the distribution of the MLEs may be far from normal, exhibiting bias, skewness, and also kurtosis that are not compatible with the normal distribution. Other standard methods for constructing confidence intervals are based on the score statistic and the signed likelihood ratio statistic, but both rely on first-order asymptotics.

A vast literature on bias adjustments for MLEs is available; see for instance \citet{FIRTH}, \citet{FERRARICRIBARI98}, \citet{CORDEIROCRIBARI}, \citet{PAGUI}, and references therein. For a review paper, see \citet{KOSMIDIS}.
Additionally, effort has been made to correct the coverage of approximate confidence intervals in limited samples. \citet{BARTLETT} proposed an approximate confidence interval based on a skewness correction to the score function. An alternative possible approach is the use of computer intensive methods such as the bootstrap procedure \citep{DICICCIOEFRON}. 

This paper proposes accurate approximate confidence intervals for a scalar parameter of interest possibly in the presence of a vector of nuisance parameters in general parametric families. We construct the proposed confidence intervals from a third-order approximation to the quantiles of the score function. Some 
features of the proposed confidence intervals are: 
first, the coverage approximation error is $O(n^{-3/2})$ while the approximation error of confidence intervals based on the first-order asymptotic distribution of the Wald, score, and signed likelihood ratio statistics is $O(n^{-1/2})$; 
second, they are equivariant under interest-respecting reparameterizations; 
third, they account for skewness and kurtosis of the distribution of the score function; 
fourth, they do not require computer-intensive resampling mechanisms; 
fifth, the confidence limits are simply computed from modified score equations. For a review on different routes for achieving third-order accuracy in confidence intervals for a scalar parameter of interest, including small-sample likelihood asymptotics, bootstrapping and the objective Bayes approach, the reader is referred to \citet{YOUNG2009}.

In Section \ref{sec:nonuisance}, we deal with the one parameter case. In Section \ref{sec:nuisance}, we extend the results to the case of a single parameter of interest and a fixed number of nuisance parameters. An implementation of accurate approximate confidence intervals for regression models is presented in Section \ref{sec:generalregression}. Monte Carlo evidence on the performance of the modified confidence intervals is presented in Section \ref{sec:MonteCarlo}. Applications are presented and discussed in Section \ref{sec:applications}. The paper ends with concluding remarks in Section~\ref{sec:conclusion}. Some technical details are left for three appendices.

\section{One parameter}\label{sec:nonuisance}

Let $y$ be the data, for which we consider a parametric model indexed by a parameter $\theta \in \Theta \subset \mathbb{R}$. 
Let $\ell(\theta)$ be the log-likelihood function and $U(\theta)=\partial \ell(\theta) /\partial \theta$ be the score function.
The first four cumulants of $U(\theta)$ are 
$\kappa_1(\theta)={\rm E}_\theta(U(\theta)),$
$\kappa_2(\theta)={\rm E}_\theta(U(\theta)^2),$
$\kappa_3(\theta)={\rm E}_\theta(U(\theta)^3),$ and
$\kappa_4(\theta)={\rm E}_\theta(U(\theta)^4)-3\kappa_2^2.$
We assume the usual conditions for regular models as described in \citet[Sect. 3.4]{SEVERINI2000}.
For regular models the cumulants above are finite and they are of order $n$, where $n$ is the sample size. 
Also, the first four Bartlett identities hold; in particular 
$\kappa_1(\theta)=0$ and $\kappa_2(\theta)={\rm var}_\theta(U(\theta))=i(\theta)$, where
$i(\theta)$ is the Fisher information \citep[Sect. 3.5]{SEVERINI2000}.
For simplicity of notation, from now on we drop the argument $\theta$ and denote the cumulants by $\kappa_1,\ldots,\kappa_4$.

The maximum likelihood estimator of $\theta$, $\widehat\theta$, is the solution of the estimating equation $U(\theta)=0$, and is approximately normally distributed with mean $\theta$ and variance $\kappa_2^{-1}=i(\theta)^{-1}$. Let $u_\alpha$ be the $\alpha$-quantile of a standard normal distribution
and let $\widehat\theta_\alpha=\widehat\theta-u_\alpha/\surd{i(\widehat\theta)}$. For $\alpha<0.5$, $(-\infty, \; \widehat\theta_\alpha]$ and $[\widehat\theta_{1-\alpha},\; +\infty)$ are the usual Wald-type one-sided confidence intervals for $\theta$ with approximate confidence level $\gamma=1-\alpha$. Analogously, $[\widehat\theta_{1-\alpha/2},\; \widehat\theta_{\alpha/2}]$ is  the usual Wald-type two-sided confidence interval for $\theta$ with approximate confidence level $\gamma=1-\alpha$. These confidence intervals are first-order accurate; the approximation error is $O(n^{-1/2})$. 

The Cornish-Fisher expansion of the $\alpha$-quantile of the score function is
\begin{equation*}\label{expansionqU}
q_\alpha(U(\theta)) = u_\alpha \sqrt{\kappa_2}  
+\frac{1}{6}\frac{\kappa_3}{\kappa_2}(u_\alpha^2-1)
+\frac{1}{24}\frac{\kappa_4}{\kappa_2^{3/2}}(u_\alpha^3-3u_\alpha) 
-\frac{1}{36}\frac{\kappa_3^2}{\kappa_2^{5/2}}(2u_\alpha^3-5u_\alpha)
+O(n^{-1})
\end{equation*}
\citep[Chap. 10, eq. (10.19)]{PACE}.
We then define the $\alpha$-quantile modified score 
\begin{equation}
\label{modifiedscore}
\widetilde U_\alpha(\theta)=
U(\theta)
-u_\alpha \sqrt{\kappa_2}
-\frac{1}{6}\frac{\kappa_3}{\kappa_2}(u_\alpha^2-1)
-\frac{1}{24}\frac{\kappa_4}{\kappa_2^{3/2}}(u_\alpha^3-3u_\alpha) 
+\frac{1}{36}\frac{\kappa_3^2}{\kappa_2^{5/2}}(2u_\alpha^3-5u_\alpha).
\end{equation}
We have that $q_\alpha(\widetilde U_\alpha(\theta))=O(n^{-1})$ and 
\begin{equation}
\label{modifiedscorefda}
P_\theta(\widetilde U_\alpha(\theta) \leq 0) = \alpha + O(n^{-3/2});
\end{equation}
see Appendix~\ref{ape:approximationorder}, for a proof.

If $\widetilde{\theta}_\alpha$ is the unique solution of the $\alpha$-quantile modified score equation 
\begin{equation}
\label{modifiedscoreeq}
\widetilde{U}_\alpha(\theta)=0,
\end{equation}
and $\widetilde U_\alpha(\theta)>0$ for $\theta<\widetilde\theta_\alpha$, then the events $\widetilde U_\alpha(\theta)\leq 0$ and $\widetilde \theta_\alpha\leq \theta$ are equivalent and
\begin{equation}
\label{estimatorfda}
P_\theta(\widetilde \theta_\alpha \leq \theta ) = \alpha + O(n^{-3/2}).
\end{equation}
Hence, $\theta$ is approximately the $\alpha$-quantile of the distribution of $\widetilde \theta_\alpha$; the approximation error is of order $n^{-3/2}$.

\paragraph{\bf Remark 1}. The behavior of the modified score function $\widetilde U_\alpha(\theta)$ is mostly governed by the leading term in (\ref{modifiedscore}), i.e. the score function, because the remaining terms are of smaller order. If the likelihood function is log-concave with a single maximum, the MLE ($\widehat\theta$) is the unique solution of the score equation, $U(\theta)=0$, and $U(\theta)>0$ 
for $\theta<\widehat\theta$. It is clear, however, that for a fixed small sample size, combined with $\alpha$ close to zero or one, the behavior of the observed $\widetilde U_\alpha(\theta)$ may be impacted by the remaining terms. In such a case, the $\alpha$-quantile score equation may have no solution or more than one solution. In this connection, see Example 3 below. In the Monte Carlo experiments reported in Section \ref{sec:MonteCarlo} with thousands of limited sized samples and $\alpha$ ranging from $0.5\%$ to $99.5\%$,  standard root-finding algorithms typically converge (in the worst case, the algorithm employed achieved convergence for $99.97\%$ of the simulated samples).
\vspace{0.3cm}

\paragraph{\bf Remark 2.} When $\alpha=0.5$, we have $u_{0.5}=0$ and (\ref{modifiedscore}) and (\ref{estimatorfda}) reduce to
$$\widetilde U_{0.5}(\theta) = U(\theta) + \frac{1}{6}\frac{\kappa_3}{\kappa_2},$$
and $P_\theta(\widetilde \theta_{0.5} \leq \theta ) = 0.5 + O(n^{-3/2})$,
which coincide with eq. (1) and (3) in \citet{PAGUI}, respectively. Hence, $\widetilde \theta_{0.5}$ is a third-order median bias reduced estimator of $\theta$ and it is the same as that derived in \citet{PAGUI}.
\vspace{0.3cm}

Equation (\ref{estimatorfda}) can be used to obtain confidence limits for $\theta$. For $\alpha<0.5$, we have that $(-\infty, \widetilde \theta_\alpha]$ and $[\widetilde \theta_{1-\alpha},+\infty)$ are one-sided confidence intervals for $\theta$ with approximate confidence level $\gamma=1-\alpha$. Analogously, $[\widetilde \theta_{1-\alpha/2},\widetilde \theta_{\alpha/2}]$ is a two-sided confidence interval for $\theta$ with approximate confidence level $\gamma=1-\alpha$. These confidence intervals are third-order accurate, i.e the approximation absolute error is $O(n^{-3/2})$, This is an improvement over the $O(n^{-1/2})$ approximation error of confidence intervals based on the first-order asymptotic distribution of the Wald, score, and signed likelihood ratio statistics.  We shall name the confidence intervals proposed here quantile bias reduced (QBR) confidence intervals.

\paragraph{\bf Remark 3.} $\widetilde \theta_\alpha$ is equivariant under reparameterization. Let $\omega(\theta)$ be a smooth reparameterization with inverse $\theta(\omega)$. In the parameterization $\omega$, the score function is $U(\theta(\omega)) \theta'(\omega)$, and its first four cumulants are  $\kappa_r \theta'(\omega)^r$, $r=1, 2, 3, 4,$ where $\theta'(\omega)$ is the derivative of $\theta(\omega)$ with respect to $\omega$.
Hence, the $\alpha$-quantile modified score in the parameterization $\omega$ is  $\widetilde U_\alpha(\theta(\omega))  \theta'(\omega)$, and the solution of $\widetilde U_\alpha(\theta(\omega))  \theta'(\omega) = 0$ is $\widetilde\omega_\alpha=\omega(\widetilde\theta_\alpha)$. 
\vspace{0.3cm}

\paragraph{\bf Example 1.} {\it One parameter exponential family.}
For a random sample $y_1,\ldots,y_n$ of a one parameter exponential family with canonical parameter $\theta$ and pdf
$$
f(y;\theta)=\exp\{ \theta T(y) - A(\theta)\}h(y),
$$
we have $U(\theta)= \sum_{i=1}^{n}T(y_i) - n \; d A(\theta) / d \theta$, $\kappa_r= n \; d^r A(\theta) / d \theta^r$, $r=2,3,4$. The modified score equation is given as in (\ref{modifiedscoreeq}) by plugging these quantities in (\ref{modifiedscore}).

For a random sample of an exponential distribution with mean $1/\theta$, we have 
$T(y) = -y$, $A(\theta)=-\log(\theta)$, $U(\theta)=n/\theta - \sum_{i=1}^n y_i,$ 
$\widehat\theta=n/\sum_{i=1}^n y_i$, $\kappa_2=n / \theta ^2$, $\kappa_3=-2n / \theta^3,$ $\kappa_4=6n /\theta^4$, 
$
\widehat\theta_\alpha = \widehat\theta(1-u_\alpha/\surd n),
$
and 
$
\widetilde \theta_\alpha = \widehat\theta(1-c_{n,\alpha}/\surd n),
$
where $c_{n,\alpha}=u_\alpha - (u_\alpha^2-1)/(3\surd n) + (u_\alpha^3-7u_\alpha)/(36 n).$ 
A third-order median bias reduced estimator of $\theta$ is $\widetilde \theta_{0.5} = \widehat\theta(1-1/(3n))$. 

\begin{figure}[!ht]
	\centering
		\includegraphics[width=5.55cm,height=5.55cm]{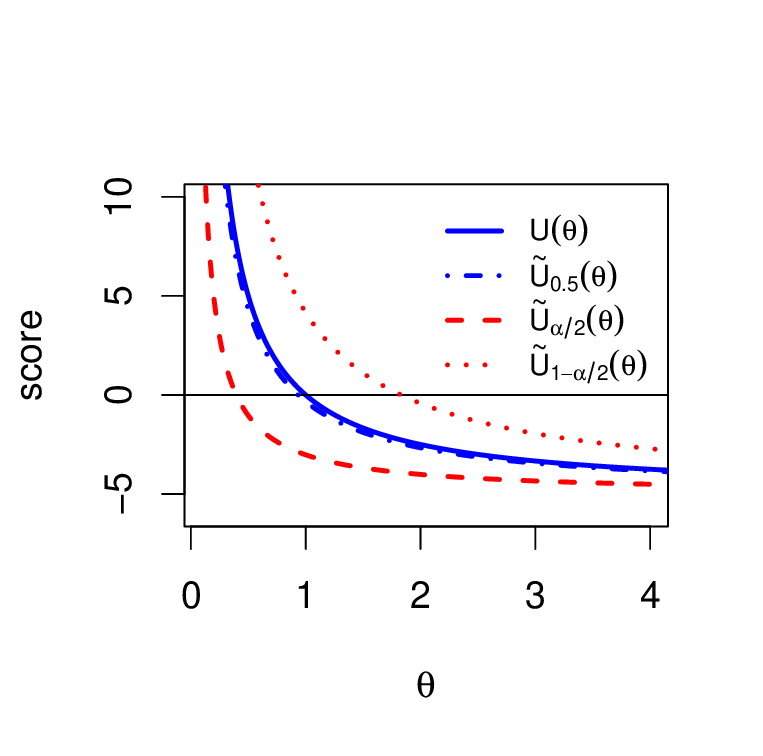}
		\includegraphics[width=5.55cm,height=5.55cm]{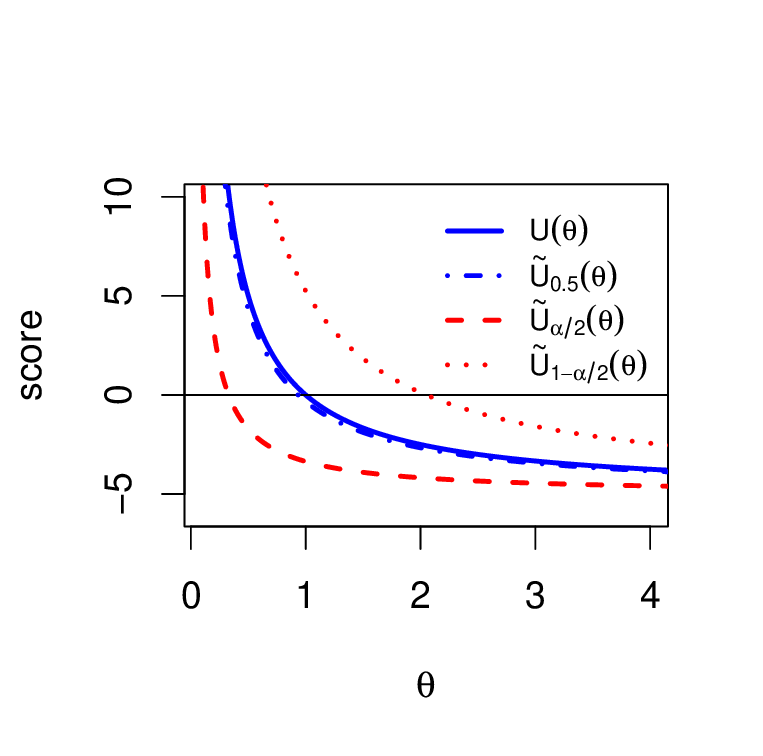}
		\includegraphics[width=5.55cm,height=5.55cm]{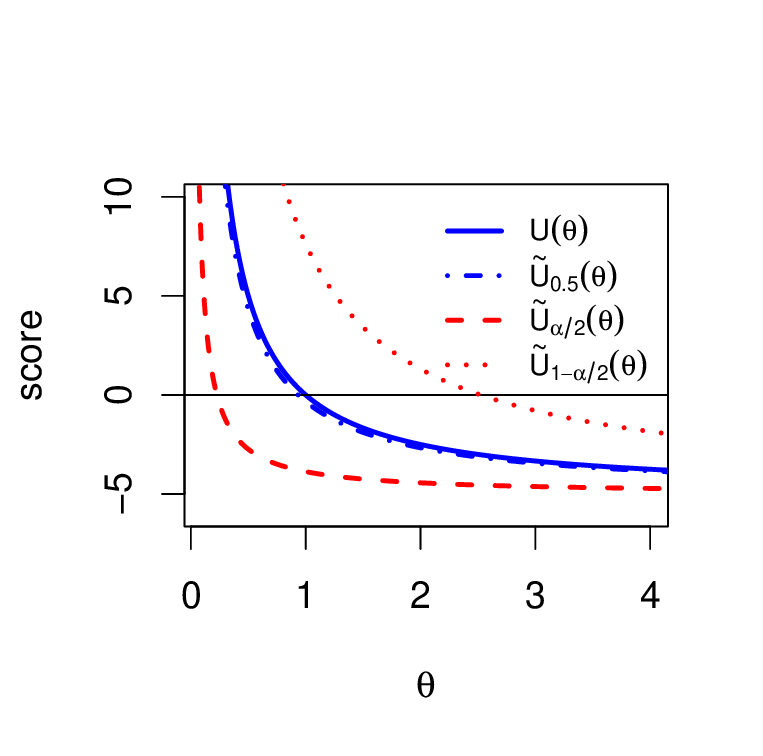}
	\caption{\footnotesize Plots of the score and modified score functions with $1-\alpha=90\%$ (left), $95\%$ (middle), and $99\%$ (right); exponential distribution, $n=5$. } 
\label{fig:scoreExp}
\end{figure}

Table \ref{tab:exponential} presents approximate 90\%, 95\%, and 99\% confidence intervals for $\theta$ based on samples with mean equal to 1, i.e. the MLE of $\theta$ is $\widehat\theta=1$, for different small sample sizes, namely $n=3, 5, 7.$ The approximate confidence limits considered are the Wald-type confidence limits (ML), 
an adjusted version that replaces $\widehat\theta$ by the median bias reduced estimator $\widetilde\theta_{0.5}$ (MBR),  the usual score confidence limits that come from the first two terms in the right hand side of (\ref{modifiedscore}) (score), and the confidence limits proposed here (QBR). In this case, the ML and score confidence limits coincide.
Like the ML and score confidence limits, the MBR confidence limits are first-order accurate, while those proposed in this paper are third-order accurate. For comparison, the table includes confidence intervals derived from the signed log likelihood ratio statistic ($r$), and its third-order modification ($r^*$); \citet{BARNDORFF1986}. The table also gives the exact $1-\alpha$ confidence intervals $[\chi^2_{\alpha/2,2n}/(2\sum_{i=1}^n y_i), \chi^2_{1-\alpha/2,2n}/(2\sum_{i=1}^n y_i)]$, where $\chi^2_{\alpha,n}$ is the $\alpha$-quantile of a chi-squared distribution with $n$ degrees of freedom. The figures in Table \ref{tab:exponential} reveal that QBR confidence limits proposed here and those derived from $r^*$ are remarkably accurate even in very small samples. They lead to a considerable improvement over both first-order approximate confidence intervals. The ML, score, and MBR confidence intervals may even include negative numbers, i.e. values outside the parameter space. Figure \ref{fig:scoreExp} shows plots of the score and modified score functions for illustrative purposes.  

\vspace{0.3cm}
\begin{table}[!ht]
\begin{center}
\caption{\footnotesize Approximate and exact $1-\alpha$ confidence intervals for $\theta$ for different sample sizes $n$; exponential distribution} \label{tab:exponential}
{\footnotesize
\begin{tabular}{llrrr}  \hline
                             &                     &\multicolumn{3}{c}{$1-\alpha$}\\               
\cline{3-5}
 $n$                         &                     & {90\%}           & {95\%}            & {99\%}\\ \hline
\multirow{3}{*}{\centering 3}&  ML / score         & $[0.05, \ 1.95]$ & $[-0.13, \ 2.13]$ & $[-0.49, \ 2.49]$\\
                             &  MBR                & $[0.04, \ 1.73]$ & $[-0.12, \ 1.89]$ & $[-0.43, \ 2.21]$\\
                             &  QBR                & $[0.28, \ 2.10]$ & $[0.22, \ 2.41]$  & $[0.14, \ 3.11]$\\
                             &  $r$                & $[0.32, \ 2.27]$ & $[0.25, \ 2.59]$  & $[0.14, \ 3.30]$\\
                             &  $r^*$              & $[0.27, \ 2.10]$ & $[0.21, \ 2.41]$  & $[0.11, \ 3.09]$\\
														 & Exact               & $[0.27, \ 2.10]$ & $[0.21, \ 2.41]$  & $[0.11, \ 3.09]$\\
\\															
\multirow{3}{*}{\centering 5}& ML / score          & $[0.26, \ 1.74]$ & $[0.12, \ 1.88]$  & $[-0.15, \ 2.15]$\\
                             & MBR                 & $[0.25, \ 1.62]$ & $[0.12, \ 1.75]$  & $[-0.14, \ 2.01]$\\
                             & QBR                 & $[0.40, \ 1.83]$ & $[0.33, \ 2.05]$  & $[0.23, \ 2.53]$\\
                             & $r$                 & $[0.43, \ 1.93]$ & $[0.36, \ 2.15]$  & $[0.24, \ 2.63]$\\
                             & $r^*$               & $[0.39, \ 1.83]$ & $[0.32, \ 2.05]$  & $[0.22, \ 2.52]$\\
														 & Exact               & $[0.39, \ 1.83]$ & $[0.32, \ 2.05]$  & $[0.22, \ 2.52]$\\
\\														
\multirow{3}{*}{\centering 7}& ML / score          & $[0.38, \ 1.62]$ & $[0.26, \ 1.74]$  & $[0.03, \ 1.97]$\\
                             & MBR                 & $[0.36, \ 1.54]$ & $[0.25, \ 1.66]$  & $[0.03, \ 1.88]$\\
                             & QBR                 & $[0.47, \ 1.69]$ & $[0.40, \ 1.87]$  & $[0.30, \ 2.24]$ \\
                             & $r$                 & $[0.50, \ 1.76]$ & $[0.43, \ 1.93]$  & $[0.31, \ 2.31]$\\
                             & $r^*$               & $[0.47, \ 1.69]$ & $[0.40, \ 1.87]$  & $[0.29, \ 2.24]$\\
														 & Exact               & $[0.47, \ 1.69]$ & $[0.40, \ 1.87]$  & $[0.29, \ 2.24]$\\
\hline
\end{tabular}
}
\end{center} 
\end{table}	
																									
\paragraph{\bf Example 2.} {\it Normal variance.}
For a random sample $y_1,\ldots,y_n$ of a normal distribution with mean zero and variance $\theta$, the pdf for a single observation is
$
f(y;\theta)=(2\pi\theta)^{-1/2}\exp\{ -{y^2}/(2\theta)\}, 
$
$y \in \mathbb{R}$, 
and we have $U(\theta)= n/(2\theta)(\widehat\theta/\theta-1)$, where $\widehat\theta=\sum_{i=1}^{n}y_i^2/n$, $\kappa_2= n/(2\theta^2)$, $\kappa_3= n/\theta^3$, and $\kappa_4= 3n/\theta^4$. The modified score equation is given as in (\ref{modifiedscoreeq}) by plugging these quantities in (\ref{modifiedscore}). It leads to
$
\widetilde U_\alpha(\theta) = n/(2\theta)(\widehat\theta/\theta-1-k_{n,\alpha}/\sqrt{n}),
$ 
where 
$
k_{n,\alpha}=\surd{2}[u_\alpha +\surd{2} (u_\alpha^2-1)/(3\surd n) + (u_\alpha^3-7u_\alpha)/(18 n)].
$
Then $\widetilde\theta_\alpha=\widehat\theta/(1+k_{n,\alpha}/\surd n)$, that is the unique root of the modified score equation,
and $\widetilde U_\alpha(\theta)>0$ for  $\theta<\widetilde\theta_\alpha$, if $1+k_{n,\alpha}/\surd n > 0$. 
A third-order median bias reduced estimator of $\theta$ is 
$
\widetilde \theta_{0.5} = \widehat\theta/(1-2/(3 n)).
$ 

\begin{figure}[!ht]
	\centering
		\includegraphics[width=5.55cm,height=5.55cm]{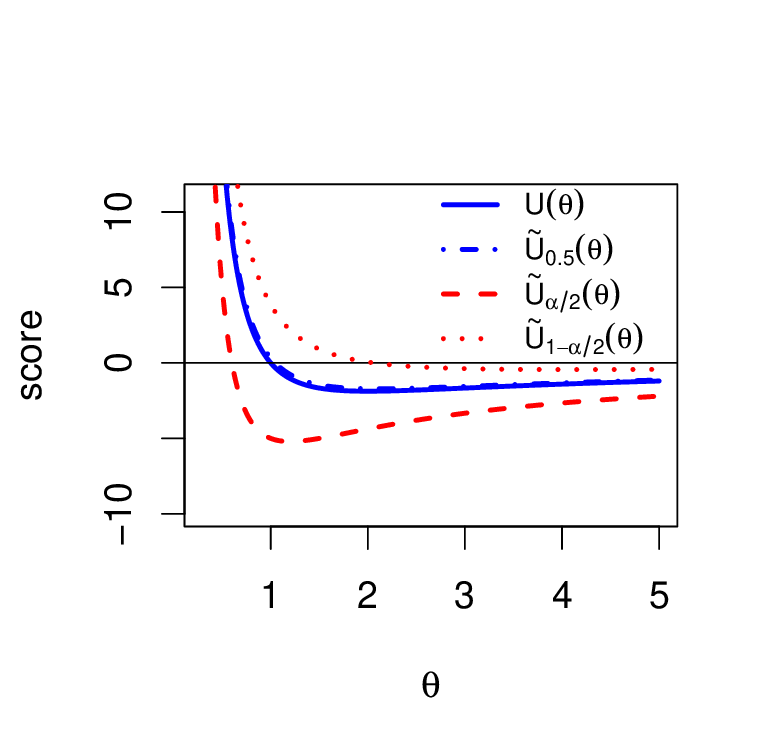}
		\includegraphics[width=5.55cm,height=5.55cm]{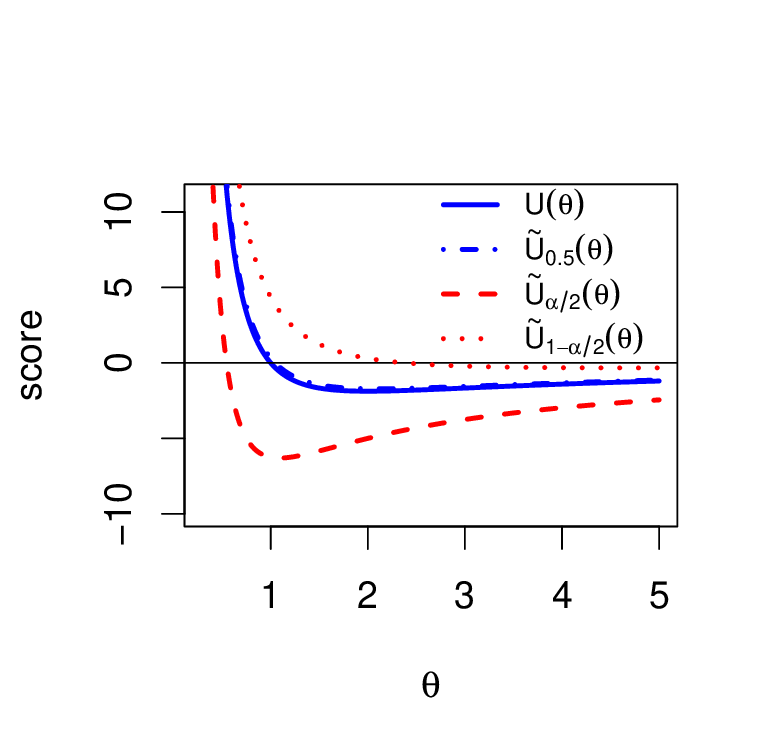}
		\includegraphics[width=5.55cm,height=5.55cm]{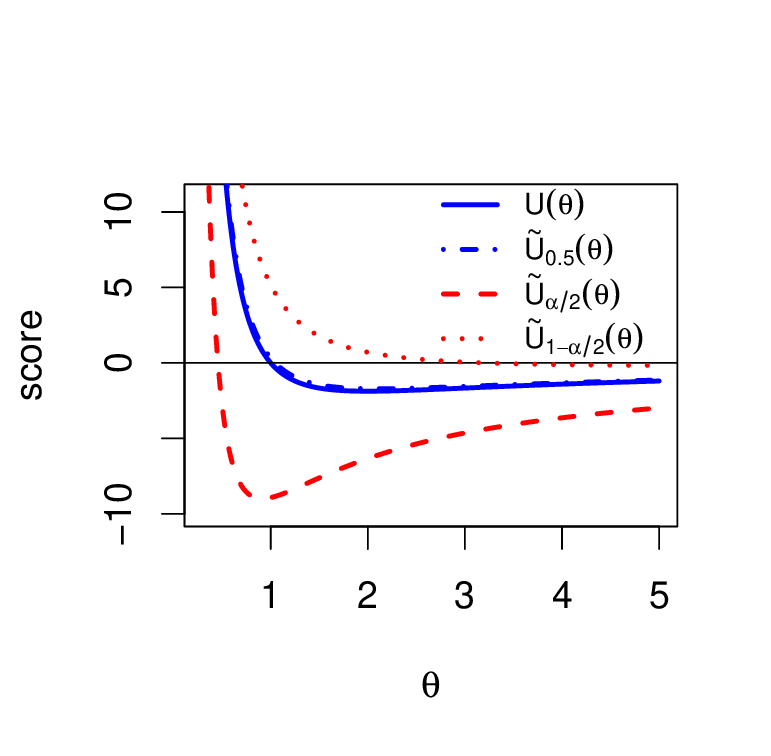}
	\caption{\footnotesize Plots of the score and modified score functions with $1-\alpha=90\%$ (left), $95\%$ (middle), and $99\%$ (right); normal variance, $n=15$.} 
\label{fig:scoreNormalVariance}
\end{figure}

Table \ref{tab:NormalVariance} presents approximate 90\%, 95\%, and 99\%  confidence intervals for $\theta$ based on samples of sizes $n=10, 15, 20$ with $\widehat\theta=1$. The exact $1-\alpha$ confidence interval is
$
[\sum_{i=1}^n y_i^2/\chi^2_{1-\alpha/2,n},$ $ \sum_{i=1}^n y_i^2/\chi^2_{\alpha/2,n}].
$  
The third-order approximate confidence limits proposed here and those derived from $r^*$ are much more accurate than the first-order approximate confidence intervals. Figure \ref{fig:scoreNormalVariance} shows plots of the score and modified score functions for the sample size $n=15$.

\begin{table}[!ht]
\begin{center}
\caption{\footnotesize Approximate and exact $1-\alpha$ confidence intervals for $\theta$ for different sample sizes $n$; normal variance} \label{tab:NormalVariance}
{\footnotesize
\begin{tabular}{llrrr}  \hline
                              &                       &\multicolumn{3}{c}{$1-\alpha$}\\               
\cline{3-5}
 $n$                          &         & {90\%}           & {95\%}            & {99\%}\\ \hline
\multirow{3}{*}{\centering 10}&  ML     & $[0.26, \ 1.74]$ & $[0.12, \ 1.88]$  & $[-0.15, \ 2.15]$\\
                              &  MBR    & $[0.28, \ 1.86]$ & $[0.13, \ 2.01]$  & $[-0.16, \ 2.31]$\\
															&  score  & $[0.58, \ 3.78]$ & $[0.53, \ 8.10]$  & $[0.46, \ \infty]$\\
                              &  QBR    & $[0.55, \ 2.53]$ & $[0.49, \ 3.05]$  & $[0.40, \ 4.42]$\\
                              &  $r$    & $[0.52, \ 2.31]$ & $[0.47, \ 2.79]$  & $[0.38, \ 4.15]$\\
                              &  $r^*$  & $[0.55, \ 2.54]$ & $[0.49, \ 3.08]$  & $[0.40, \ 4.64]$\\
														  & Exact   & $[0.55, \ 2.54]$ & $[0.49, \ 3.08]$  & $[0.40, \ 4.64]$\\
\\															
\multirow{3}{*}{\centering 15}&  ML     & $[0.40, \ 1.60]$ & $[0.28, \ 1.72]$  & $[0.06, \ 1.94]$\\
                              &  MBR    & $[0.42, \ 1.68]$ & $[0.30, \ 1.80]$  & $[0.06, \ 2.03]$\\
															&  score  & $[0.62, \ 2.50]$ & $[0.58, \ 3.52]$  & $[0.52, \ 16.82]$\\
                              &  QBR    & $[0.60, \ 2.06]$ & $[0.55, \ 2.39]$  & $[0.46, \ 3.21]$\\
                              &  $r$    & $[0.58, \ 1.95]$ & $[0.53, \ 2.25]$  & $[0.44, \ 3.05]$\\
                              &  $r^*$  & $[0.60, \ 2.07]$ & $[0.55, \ 2.40]$  & $[0.46, \ 3.26]$\\	
														  & Exact   & $[0.60, \ 2.07]$ & $[0.55, \ 2.40]$  & $[0.46, \ 3.26]$\\
\\														
\multirow{3}{*}{\centering 20}&  ML     & $[0.48, \ 1.52]$ & $[0.38, \ 1.62]$  & $[0.19, \ 1.81]$\\
                              &  MBR    & $[0.50, \ 1.57]$ & $[0.39, \ 1.68]$  & $[0.19, \ 1.88]$\\
															&  score  & $[0.66, \ 2.08]$ & $[0.63, \ 2.62]$  & $[0.55, \ 5.39]$\\
                              &  QBR    & $[0.64, \ 1.84]$ & $[0.59, \ 2.08]$  & $[0.50, \ 2.67]$ \\
                              &  $r$    & $[0.62, \ 1.77]$ & $[0.57, \ 2.00]$  & $[0.49, \ 2.56]$\\
                              &  $r^*$  & $[0.64, \ 1.84]$ & $[0.58, \ 2.09]$  & $[0.50, \ 2.69]$\\	
														  & Exact   & $[0.64, \ 1.84]$ & $[0.59, \ 2.09]$  & $[0.50, \ 2.69]$\\
\hline
\end{tabular}
} 
\end{center}         
\end{table}	
																									
\paragraph{\bf Example 3.} {\it Skew-normal distribution.} Let $y_1,\ldots,y_n$ be independent random variables from a skew-normal distribution with shape parameter $\theta \in  \mathbb{R}$ and pdf $f(y;\theta)= 2\varphi(y)\Phi(\theta y)$, $y \in \mathbb{R}$, where $\varphi(\cdot)$ and $\Phi(\cdot)$ denote the standard normal density and distribution functions, respectively. The score function is $U(\theta) = \sum_{i=1}^{n}\zeta(\theta y_i)y_i$, where $\zeta(x)={\rm d} \log\{2\Phi(x)\}/{\rm d} x$. Let $a_{kl}={\rm E}_\theta(y^k \zeta(\theta y)^l)$. The expected quantities needed to compute the modified score (\ref{modifiedscore}) are $\kappa_2=na_{22}$, $\kappa_3=na_{33}$, and $\kappa_4=n(a_{44}-3a_{22}^2)$.

\begin{figure}[!ht]
	\centering
		\includegraphics[width=6cm,height=6cm]{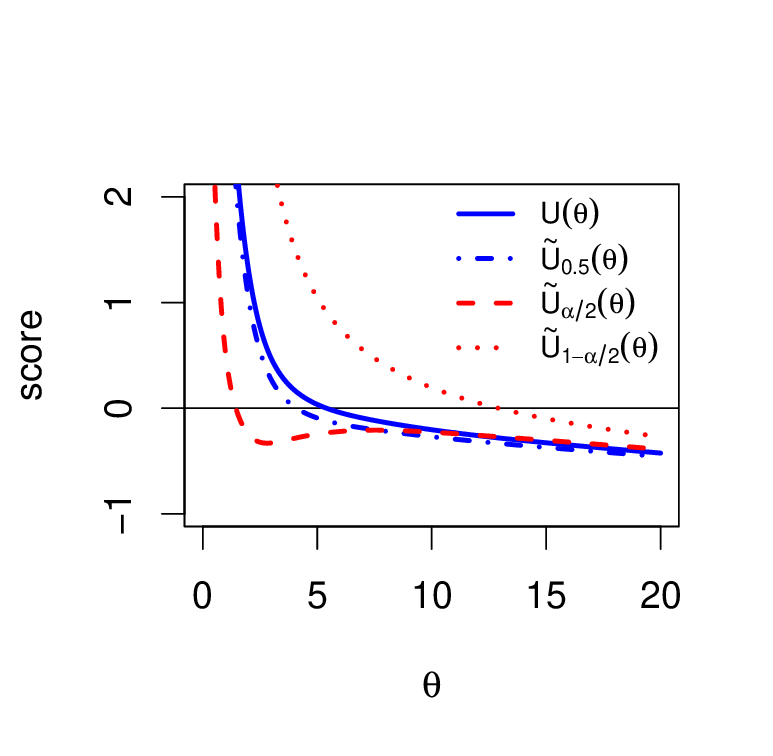} \hspace{2cm}
		\includegraphics[width=6cm,height=6cm]{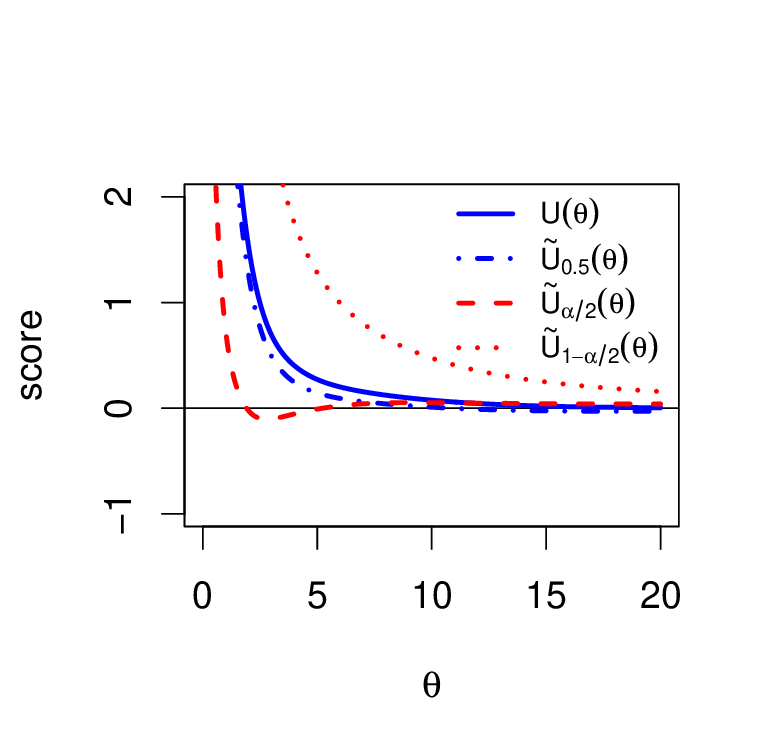}
	\caption{\footnotesize Plots of the score and modified score functions for the data in Example 1 of \citet{SARTORI} with $1-\alpha=95\%$; original data (left), modified data (right).}
\label{fig:scoreSNSartori}
\end{figure}

The maximum likelihood estimator may deliver infinite estimates with positive probability. This is illustrated in Example 1 of \citet{SARTORI} (see also Figure 1 in \citet{PAGUI}), in which the change of sign of the only negative observation moves the MLE from $5.40$ to $\infty$. Figure \ref{fig:scoreSNSartori} shows the plots of the score and modified score functions for the original and modified data with $1-\alpha=95\%$. For the original sample, the  ML and MBR confidence intervals are $[-0.64; 11.44]$ and $[0.02; 8.17]$ for the confidence level $95\%$, while the QBR confidence interval proposed here is $[1.44;  12.91]$. For the modified sample, the ML confidence interval cannot be obtained because the MLE is infinite. The MBR confidence interval is $[-5.93; 27.58]$, that is uninformative because $\theta$ in this range covers from highly negative skewed to highly positive skewed normal distributions. The modified score function $\widetilde U_{0.025}(\theta)$ crosses the horizontal axis twice. The shape of the $\widetilde U_{0.025}(\theta)$ plot suggests the choice of the smallest root in this case. Moreover, $\widetilde U_{0.025}(\theta)$ crosses the horizontal axis from positive to negative values only at this root. By taking the smallest root, the confidence interval is $[1.91; \infty)$, suggesting a highly positive skew normal distribution. 
In this extreme case, for which the MLE is infinite and the  ML confidence interval cannot be obtained, the modified score function leads to reasonable inference.

As problems in the estimation of $\theta$ tend to arise when $|\theta|$ is far from zero, a reviewer suggested to perform a simulation study using the settings of Example 2 of \citet{PAGUI}, i.e. with $\theta=5, 10$, $n=20, 50, 100$, and the $95\%$ nominal coverage probability, and to report how often it is not possible to compute the confidence interval proposed here. Confidence intervals that have at least one finite end-point are considered valid.\footnote{The confidence intervals based on the signed likelihood root are always valid, because even when the MLE is $\infty$ ($-\infty$) the likelihood converge to a constant value and a lower (upper) bound for a confidence interval can always be found.} Table \ref{tab:skewNormalSimul} shows the results based on 10,000 Monte Carlo simulated samples. The QBR confidence interval can be computed more often, in some instances much more often, than the ML confidence interval. For example, when $n=20$ and $\theta=10$, the MLE estimate of $\theta$ is finite in only $47.3\%$ of the simulated samples; hence in $53.7\%$ of the times the ML confidence interval cannot be computed. In the same scenario, the QBR confidence interval cannot be computed in only $10.6\%$ of the times. Table \ref{tab:skewNormalSimul} presents coverage estimated probabilities computed as the percentage of the samples for which the confidence interval covers the parameter value out of the samples that lead to a valid confidence interval. Hence, coverage probabilities should be judged with caution. For instance, although the coverage probability of the ML confidence interval is not too far from 95\% when $n=20$ and $\theta=10$, it is computed out of 47.3\% of the samples only. The results shown in Table \ref{tab:skewNormalSimul} suggest that the QBR  and MBR confidence intervals outperform the ML confidence interval in terms of both coverage probability and computability. An advantage of the MBR confidence interval is that it can be computed in all the simulated samples.

\begin{table}[!ht]
\begin{center}
\caption{\footnotesize Simulation results of the interval estimation of the skew-normal shape parameter; $\%$possible, percentage of the samples for which it is possible to compute the confidence interval; $\%$coverage, percentage coverage of the confidence interval out of those that are valid. }\label{tab:skewNormalSimul}
{\footnotesize	
\begin{tabular}{rrrcrcc} \hline
$\theta$  & $n$   &                         & $\%$possible & $\%$coverage\\
    \hline													                  	   	
    5     & 20    & ML                   	  &  72.7 		   & 95.1     	 \\
          &       & MBR                 	  &  100  		   & 91.4   	   \\
          &       & QBR                 		&  95.6 		   & 94.9 		   \\
\\                  										                  	   	
          & 50    & ML                  		&  95.6 		   & 96.6 		   \\
          &       & MBR                 		&  100  		   & 94.3 		   \\
          &       & QBR                 		& 99.7  		   & 94.2 		   \\
\\ 															                  	   	
          & 100   & ML                   		& 99.8  		   & 96.7 		   \\
          &       & MBR                 		& 100   		   & 95.1 		   \\
          &       & QBR                 		& 100   		   & 94.9 		   \\
\\ 														                  	   	
 10       & 20    & ML                  		& 47.3  		   & 91.0 		   \\
          &       & MBR                 		& 100   		   & 84.4 		   \\
          &       & QBR                  		& 89.4  		   & 96.2 		   \\
\\         													                  	   	
          & 50    & ML                  		& 79.7  		   & 95.3 		   \\
          &       & MBR                 		& 100   		   & 91.6 		   \\
          &       & QBR                 		& 97.8  		   & 93.0 		   \\
\\															                  	   	
          & 100   & ML                  		& 95.9  		   & 96.3 		   \\
          &       & MBR                 		& 100   		   & 94.1 		   \\
          &       & QBR                 		& 99.9  		   & 94.2 		   \\
\hline  
    \end{tabular}%
}
\end{center}   		
\end{table}%

\section{Scalar parameter of interest and nuisance parameters}\label{sec:nuisance}

Let $\theta=(\psi,\lambda^\top)^\top$, where $\psi$ is the scalar parameter of interest and $\lambda=(\lambda_1,\ldots,\lambda_p)^\top$ is the nuisance parameter. Let the subscript $\psi$ refer to the parameter $\psi$ and the indices $a, b, c,\ldots$ refer to the components of $\lambda$, so that the log-likelihood derivatives are $U_\psi=U_\psi(\psi,\lambda)=\partial \ell(\psi,\lambda)/\partial\psi$, $U_a=U_a(\psi,\lambda)=\partial\ell(\psi,\lambda)/\partial\lambda_a$,
$
U_{ab}=\partial^2 \ell(\psi,\lambda)/\partial \lambda_a \partial \lambda_b, 
U_{\psi a}=\partial^2 \ell(\psi,\lambda)/\partial \psi \partial \lambda_a,
$
etc., $a=1,\ldots,p$. Consider the index notation for the joint cumulants of derivatives of the log-likelihood function \citep{LAWLEY,HAYAKAWA,MCCULLAGH84,MCCULLAGH87}: 
$
\kappa_\psi={\rm E}(U_\psi)=0, 
\kappa_a={\rm E}(U_a)=0,
\kappa_{ab}={\rm E}(U_{ab}),
\kappa_{abc}={\rm E}(U_{abc}), 
\kappa_{abcd}={\rm E}(U_{abcd}), 
\kappa_{a,b}={\rm E}(U_a U_b),
\kappa_{a,bc}={\rm E}(U_a U_{bc}), 
\kappa_{ab,cd}={\rm E}(U_{ab} U_{cd})-\kappa_{ab}\kappa_{cd}, 
\kappa_{a,b,cd}={\rm E}(U_a U_b U_{cd})- \kappa_{a,b}\kappa_{cd},
\kappa_{a,b,c,d}={\rm E}(U_a U_b U_c U_d) - \kappa_{a,b}\kappa_{c,d} - \kappa_{a,c}\kappa_{b,d} - \kappa_{a,d}\kappa_{b,c},
\kappa_{\psi,ab}={\rm E}(U_\psi U_{ab}),
$
etc. All $\kappa$'s refer to a total over the sample and are of order $n$.
In addition, $\kappa^{a,b}$ denotes the element $(a,b)$ of the inverse covariance matrix of $U_\lambda$. In the sequel, the Einstein convention of sum of repeated indices is used, i.e. if an index occurs both as a superscript and as a subscript in a single term then summation over that index is understood.  
 
Let $\widehat\theta=(\widehat\psi,\widehat\lambda^\top)^\top$ be the MLE of $\theta$ and let $\widehat\lambda_\psi$ be the MLE of $\lambda$ for fixed $\psi$.
 The profile log-likelihood function for $\psi$ is
$$
\ell_P(\psi)=\ell(\psi,\widehat\lambda_\psi^\top),
$$
and the score function derived from the profile log-likelihood function is
$$
U_P(\psi)=\frac{\partial \ell_P(\psi)}{\partial \psi}=U_\psi(\psi,\widehat\lambda_\psi^\top).
$$
The leading term of the expansion of $U_P(\psi)$ is the efficient score function for $\psi$, i.e. $\overline U_\psi=U_\psi - \beta_\psi^a U_a$, where $\beta_\psi^a=\kappa_{\psi,b} \kappa^{a,b}$  \citep[Example 6]{PACE2004}. 
Approximate expressions for the first four cumulants of $\overline U_\psi$  are
\begin{eqnarray}\label{profilecumulants}
\kappa_{1\psi}&=&-\frac{1}{2}\kappa^{a,b}  
\left[
\kappa_{\psi,ab} + \kappa_{\psi,a,b}
- \beta_\psi^c
\left(\kappa_{c,ab} + \kappa_{c,a,b} \right) 
\right],\nonumber\\
\kappa_{2\psi}&=& \kappa_{\psi,\psi} - \beta_\psi^a \kappa_{\psi,a}, \nonumber\\
\kappa_{3\psi}&=& \kappa_{\psi,\psi,\psi} 
- 3 \beta_\psi^a \kappa_{a,\psi,\psi} 
+ 3  \beta_\psi^a \beta_\psi^b \kappa_{a,b,\psi} 
- \beta_\psi^a \beta_\psi^b \beta_\psi^c \kappa_{a,b,c}, \nonumber\\
\kappa_{4\psi}&=& 
\kappa_{\psi,\psi,\psi,\psi}
- 4 \beta_\psi^a \kappa_{a,\psi,\psi,\psi} 
+ 6 \beta_\psi^a \beta_\psi^b \kappa_{a,b,\psi,\psi} 
- 4 \beta_\psi^a \beta_\psi^b \beta_\psi^c \kappa_{a,b,c,\psi} 
+ \beta_\psi^a \beta_\psi^b \beta_\psi^c \beta_\psi^d \kappa_{a,b,c,d}; 
\end{eqnarray}
see  
\citet[eq. (7)]{PAGUI} for $\kappa_{1\psi}$, $\kappa_{2\psi}$, and $\kappa_{3\psi}$, and 
\citet[Chap. 7]{BARNDORFF31} for $\kappa_{4\psi}$. 
If $\psi$ and $\lambda$ are orthogonal parameters, we have $\kappa_{\psi,a}=0$ and $\beta_{\psi}^{a}=0$, and hence the equations in (\ref{profilecumulants}) reduce to
\begin{equation}\label{kappasortog}
\kappa_{1\psi}=-\frac{1}{2} \kappa^{a,b}(\kappa_{\psi,ab} + \kappa_{\psi,a,b}), \ \
\kappa_{2\psi}= \kappa_{\psi,\psi}, \ \
\kappa_{3\psi}= \kappa_{\psi,\psi,\psi}, \ \
\kappa_{4\psi}= \kappa_{\psi,\psi,\psi,\psi}.
\end{equation} 

\sloppy

The Cornish-Fisher expansion of the $\alpha$-quantile of $ \overline U_\psi$ is
\begin{equation*}
q_\alpha( \overline U_\psi)= 
\kappa_{1\psi} 
+ u_\alpha \sqrt{\kappa_{2\psi}}
+\frac{1}{6}\frac{\kappa_{3\psi}}{\kappa_{2\psi}}(u_\alpha^2-1)
+\frac{1}{24}\frac{\kappa_{4\psi}}{\kappa_{2\psi}^{3/2} }(u_\alpha^3-3u_\alpha)
-\frac{1}{36}\frac{\kappa_{3\psi}^2 }{\kappa_{2\psi}^{5/2}}(2u_\alpha^3-5u_\alpha)+O(n^{-1});
\end{equation*}
see \citet[eq. (10.19)]{PACE}. Since $\overline U_\psi$ is the leading term of the asymptotic expansion of $U_P(\psi)$, we define the $\alpha$-quantile modified profile score
\begin{eqnarray}\label{modifiedprofilescore}
\widetilde U_\alpha(\psi) &=& U_P(\psi) + M_{\psi,\alpha},
\end{eqnarray}
where 
\begin{eqnarray}\label{M}
M_{\psi,\alpha} = - \kappa_{1\psi} 
- u_\alpha \sqrt{\kappa_{2\psi}} 
-\frac{1}{6}\frac{\kappa_{3\psi}}{\kappa_{2\psi}}(u_\alpha^2-1)
-\frac{1}{24}\frac{\kappa_{4\psi}}{\kappa_{2\psi}^{3/2}}(u_\alpha^3-3u_\alpha)
+\frac{1}{36}\frac{\kappa_{3\psi}^2}{\kappa_{2\psi}^{5/2}}(2u_\alpha^3-5u_\alpha).
\end{eqnarray}
If $\widetilde{\psi}_\alpha$ is the unique solution of the $\alpha$-quantile modified profile score equation, 
$$\widetilde U_\alpha(\psi)=0$$ 
with $\lambda$ replaced by $\widehat\lambda_\psi$, and $\widetilde U_\alpha(\psi)>0$ for $\psi<\widetilde\psi_\alpha$,
the events $\widetilde U_\alpha(\psi)\leq 0$ and $\widetilde \psi_\alpha \leq \psi$ are equivalent so that
\begin{equation}
\label{estimatorfdanuisance}
P_\theta(\widetilde \psi_\alpha \leq \psi ) = \alpha + O(n^{-3/2}).
\end{equation}
\noindent Similarly to the one-parameter case,  $\widetilde \psi_\alpha$ can be used to obtain approximate confidence limits for $\psi$ with approximation error $O(n^{-3/2})$ in place of the $O(n^{-1/2})$ approximation error of the confidence limits obtained from the asymptotic normality of MLEs. These confidence intervals will be referred to as quantile bias reduced (QBR) confidence intervals as in the one-parameter setup.

To numerically obtain $\widetilde \psi_\alpha$, one can solve the system of estimating equations consisting of the modified score equation for $\psi$, $U_\psi(\psi,\lambda^\top)+M_{\psi,\alpha}=0$, and the score equations $U_a(\psi,\lambda^\top)=0$, for $a=1,\ldots,p$, for the components of the nuisance parameter $\lambda$.
 
\paragraph{\bf Remark 4.} As in the one-parameter setting, there is no general guarantee that the $\alpha$-quantile modified score equation has a (unique) solution and that $\widetilde U_\alpha(\psi)>0$ for $\psi<\widetilde\psi_\alpha$. These assumptions may fail for extreme values of $\alpha$ and finite sample sizes.Example 3 suggests that a visual inspection of the (profile) modified score function may be helpful to handle non-standard situations.

\paragraph{\bf Remark 5.} When $\alpha=0.5$, (\ref{modifiedprofilescore}) and (\ref{estimatorfdanuisance}) reduce to
$$\widetilde U_{0.5}(\psi) = U_P(\theta) - \kappa_{1\psi} + \frac{1}{6}\frac{\kappa_{3\psi}}{\kappa_{2\psi}},$$
which coincides with eq. (8) in \citet{PAGUI}. As in the one-parameter case, $\widetilde \psi_{0.5}$ is therefore a third-order median bias reduced estimator of $\psi$, as derived in that paper.

\paragraph{Remark 6.} 
As in the one-parameter case, $\widetilde \psi_\alpha$ is equivariant under interest-respecting reparameterization.
Let $\omega(\phi,\tau)$ be a smooth reparameterization with $\phi=\phi(\psi)$ and $\tau=\tau(\lambda)$ being one-to-one functions with inverse $\psi(\phi)$ and $\lambda(\tau)$, respectively. In the parameterization $\omega$, the score function for $\phi$ is 
$U_\psi(\psi(\phi),\lambda(\tau)) \psi^\prime(\phi)$ and the corresponding four cumulants in (\ref{profilecumulants})  are $\kappa_{\phi,ab} = \kappa_{\psi,ab} \psi^\prime(\phi)$, where $\psi^\prime(\phi)$ is the derivative of $\psi(\phi)$ with respect to $\phi$. 
Then, $M_{\phi,\alpha}=  M_{\psi,\alpha} \psi^\prime(\phi)$;
note, for instance, that $\kappa_{\phi,ab}= \kappa_{\psi,ab} \psi^\prime(\phi)$ and 
$\beta_{\phi}^{a} = \beta_\psi^a \ \psi^\prime(\phi)$.
It follows that the $\alpha$-quantile modifed score for $\phi$ is 
$U_\phi(\phi(\psi),\tau(\lambda)) + M_{\phi,\alpha}=
(U_\psi(\psi(\phi),\lambda(\tau)) + M_{\psi,\alpha})\psi^\prime(\phi)$.  Hence,
$\widetilde\phi_\alpha = \phi(\widetilde \psi_\alpha).$ 

\paragraph{\bf Example 4.}{\it Gamma distribution with mean $\mu$ and coefficient of variation $\phi^{-1/2}$, Gamma$(\mu,\phi)$.} For a random sample $y_1,\ldots,y_n$ of a Gamma$(\mu,\phi)$ distribution with pdf 
$$
f(y;\mu,\phi)=\frac{1}{\Gamma(\phi)} \left(\frac{\phi}{\mu}\right)^\phi y^{\phi-1} \exp\left( -\frac{\phi}{\mu} y \right), \quad y>0, \quad \mu>0, \quad \phi>0,
$$
the score function for $\mu$ and $\phi$ are, respectively,
$$
U_\mu(\mu,\phi)= n\frac{\phi}{\mu^2}(\overline y - \mu)
$$
and
$$
U_\phi(\mu,\phi) = - n\Psi(\phi)+ n + n\log\frac{\mu}{\phi} + \sum_{i=1}^n \log(y_i) 
         - n\frac{\overline y}{\mu},
$$
where $\Gamma(\cdot)$ is the gamma function, $\Psi(\phi)=\Gamma^\prime(\phi)/\Gamma(\phi)$, and $\overline y = \sum_{i=1}^n y_i/n$.
Here, $\mu$ and $\phi$ are orthogonal parameters and it comes from (\ref{kappasortog}) that
$\kappa_{1\mu}=0$, \ $\kappa_{2\mu}= n\phi/\mu^2,$ \ $\kappa_{3\mu}= 2n\phi/\mu^3,$ and 
$\kappa_{4\mu}= 6n\phi/\mu^4.$ Additionally, 
$\kappa_{1\phi}=n/(2\phi),$ \ $\kappa_{2\phi}=n\left(-1/\phi+\Psi^{(1)}(\phi)\right),$ \
$\kappa_{3\phi}=n\left(1/\phi^2+\Psi^{(2)}(\phi)\right),$  and $\kappa_{4\phi}=n \left(-2/\phi^3+\Psi^{(3)}(\phi)\right),$
where $\Psi^{(r)}(\phi)={\rm d}^r \Psi(\phi)/{\rm d}\phi^r$.
These cumulants may be plugged in (\ref{modifiedprofilescore})--(\ref{M}) to compute QBR confidence intervals 
for $\mu$ and $\phi$. 

\paragraph{\bf Example 5.} {\it Symmetric distributions}. Let $y_1,\ldots,y_n$ be a random sample of a 
continuous symmetric distribution S$(\mu,\phi)$ with location parameter $\mu \in \mathbb{R}$ and scale parameter $\phi>0$, and pdf 
\begin{equation*}\label{sym}
f(y;\mu,\phi) =  \frac{1}{\phi}v\left(\left(\frac{y-\mu}{\phi}\right)^2\right), \quad y\in \mathbb{R},
\end{equation*} 
for some function $v$, called the density generating function (dgf), such that $v(u)\ge 0$, for all $u\ge 0$, and
$\int_{0}^{\infty}u^{-1/2}v(u)du=1$. Some members of the symmetric class of distributions are the normal, Student-$t$, type I logistic, type II logistic, and power exponential, with corresponding dgf: 
$v(u)=(2\pi)^{-1/2}e^{-u/2}$, 
$v(u)=\nu^{\nu/2}B(1/2,\nu/2)^{-1}(\nu+u)^{-\frac{\nu+1}{2}}$, with $\nu>0$ and $B(\cdot, \cdot)$ being the beta function, 
$v(u)=ce^{-u}(1+e^{-u})^{-2}$, with $c\cong 1.4843$ being the normalizing constant, 
$v(u)=e^{-\sqrt{u}}(1+e^{-\sqrt{u}})^{-2}$, and 
$v(u)=({1}/{C(\nu)})\exp\{-\frac{1}{2}u^{1/(1+\nu)}\}$ with $-1<\nu\leq1$ and $C(\nu)=\Gamma((3+\nu)/2)2^{(3+\nu)/2}$. 
The log-likelihood function for $\theta = (\mu,\phi)^\top$ and the score functions for $\mu$ and $\phi$ are, respectively,
\[
l(\theta) = -n\log(\phi)+\sum_{i=1}^{n} \log v(\epsilon_i^2),
\]
\[
 U_{\mu}(\mu, \phi) = \phi^{-1} \sum_{i=1}^{n} w_i \epsilon_i,
\quad
 U_{\phi}(\mu, \phi) = -n\phi^{-1} + \phi^{-1} \sum_{i=1}^{n} w_i \epsilon_i^2,
\]
where  
$\epsilon_i=(y_i-\mu)/\phi$ and $w_i=-2\dd \log v(u)/\dd u|_{u=\epsilon_i^2}$. The parameters $\mu$ and $\phi$ are orthogonal and we have from (\ref{kappasortog}) that
$\kappa_{1\mu} = \kappa_{3\mu} = 0$, \ 
$\kappa_{2\mu} = n\delta_{20000}/\phi^2,$ and
$\kappa_{4\mu} = n(\delta_{40000}-3\delta^2_{20000})/\phi^4.$ 
Additionally, 
\begin{align*}\label{eq:phi-iid}
\begin{split}
\kappa_{1\phi}&= -\frac{1}{2\phi}\frac{\delta_{00101} + 2\delta_{11001}}{\delta_{20000}}, \qquad  
\kappa_{2\phi} = \frac{n}{\phi^2}(\delta_{20002}-1), \qquad 
\kappa_{3\phi}=\frac{2n}{\phi^3}(1+\delta_{11003}), \\ 
\kappa_{4\phi}&=\frac{n}{\phi^4}(\delta_{40004}+4\delta_{30003}+12\delta_{20002}-3\delta^2_{20002}-6).
\end{split}
\end{align*}
Here, $\delta_{abcde}=\mathrm{E}(s^{(1)^{a}}s^{(2)^{b}} s^{(3)^{c}}s^{(4)^{d}}z^e)$, for $a,b,c,d,e \in \{0,1,2,3,4\}$, $s^{(r)}=\dd^r s(\epsilon)/\dd \epsilon^r$ with $s(\epsilon)=\log v(\epsilon^2)$ and $\epsilon\sim S(0,1).$
For the symmetric distributions listed above the $\delta$'s are given in \cite{Uribe2007}.
These profile cumulants may be plugged in (\ref{modifiedprofilescore})--(\ref{M}) to compute QBR confidence intervals 
for $\mu$ and $\phi$.

\section{An implementation for regression models}\label{sec:generalregression}

Let $f(y; \mu, \phi)$ be the pdf of a parametric distribution with two parameters, $\mu$ and $\phi$, and let $y_1, \ldots, y_n$ be independent random variables, where each $y_i$ has pdf $f(y; \mu_i, \phi_i)$. Consider the regression model that specifies $\mu_i$ and $\phi_i$ as
\begin{equation}
\label{linkmuphi}
g(\mu_i) = x_{i}^\top\beta, \ \ h(\phi_i) = z_{i}^\top\gamma.  
\end{equation}
Here, 
$\beta = (\beta_1, \ldots, \beta_q)^{\!\top}$ and $\gamma = (\gamma_1, \ldots, \gamma_m)^{\!\top}$ are vectors of  unknown parameters ($\beta \in \mathbb{R}^q$, $\gamma \in \mathbb{R}^m$, $q+m=p<n$), and  $x_{i}=(x_{i1}, \ldots, x_{iq})^\top$ and $z_{i}=(z_{i1}, \ldots, z_{im})^{\top}$ collect observations on covariates, which are assumed fixed and known. The link functions, $g(\cdot)$ and $h(\cdot)$, are strictly monotonic and twice differentiable, and map, respectively, the range of $\mu_i$ and $\phi_i$, into $\mathbb{R}$. 
We assume that the model matrices $X$ and $Z$, with element $(i,k)$ given by $x_{ik}$ and $z_{ik}$, respectively, are of full rank, i.e. rank$(X) = q$ and rank$(Z) = m$. 

Let 
$$
U_{\mu_i}=\frac{\partial \log f(y_i;\mu_i,\phi_i)}{\partial \mu_i},
\quad
U_{\phi_i}=\frac{\partial \log f(y_i;\mu_i,\phi_i)}{\partial \phi_i},
$$
$$
U_{\mu_i\mu_i}=\frac{\partial^2 \log f(y_i;\mu_i,\phi_i)}{\partial \mu_i\partial \mu_i},
\quad 
U_{\mu_i\phi_i}=\frac{\partial^2 \log f(y_i;\mu_i,\phi_i)}{\partial \mu_i\partial \phi_i},
\quad
U_{\phi_i\phi_i}=\frac{\partial^2 \log f(y_i;\mu_i,\phi_i)}{\partial \phi_i \partial \phi_i}.
$$
Let $\ell(\theta)=\sum_{i=1}^n \log f(y_i;\mu_i,\phi_i)$ be the log-likelihood function for $\theta=(\beta^\top,\gamma^\top)^\top$.
The first and second order log-likelihood derivatives with respect  to the unknown parameters are
\begin{eqnarray*}
\label{loglikderivatives}
U_{\beta_r}
&=& \sum_{i=1}^n x_{ir} \frac{1}{g^\prime(\mu_i)} U_{\mu_i}, \quad 
U_{\gamma_r} = \sum_{i=1}^n z_{ir} \frac{1}{h^\prime(\phi_i)} U_{\phi_i}, \\
U_{\beta_s\beta_t}
&=& \sum_{i=1}^n \left\{ x_{is} x_{it} \frac{1}{g^\prime(\mu_i)^2} U_{\mu_i\mu_i} 
- x_{is} x_{it} \frac{g^{\prime\prime}(\mu_i)}{g^\prime(\mu_i)^3} U_{\mu_i}  \right\}, \\ \nonumber
U_{\beta_s\gamma_t}
&=&\sum_{i=1}^n x_{is} z_{it}\frac{1}{g^\prime(\mu_i)} \frac{1}{h^\prime(\phi_i)} U_{\mu_i\phi_i}, \\ \nonumber
U_{\gamma_s\gamma_t}
&=& \sum_{i=1}^n  \left\{z_{is} z_{it}\frac{1}{h^\prime(\phi_i)^2} U_{\phi_i\phi_i}
 - z_{is} z_{it} \frac{h^{\prime\prime}(\phi_i)}{h^\prime(\phi_i)^3} U_{\phi_i} \right\}.
\end{eqnarray*}

Cumulants of the log-likelihood derivatives may be obtained from cumulants of derivatives of $\log f(y_i;\mu_i,\phi_i)$. For instance,
$$
\kappa_{{\beta_s,\beta_t}} = {\rm E}(U_{\beta_s} U_{\beta_t}) 
= \sum_{i=1}^n \sum_{j=1}^n x_{is} x_{it} \frac{1}{g^\prime(\mu_i)^2} {\rm E}(U_{\mu_i}U_{\mu_j})  
=\sum_{i=1}^n x_{is} x_{it} \frac{1}{g^\prime(\mu_i)^2} \kappa_{\mu_i,\mu_i},
$$
where $\kappa_{\mu_i,\mu_i}={\rm E}(U_{\mu_i}^2)$, since ${\rm E}(U_{\mu_i}U_{\mu_j})=0$, for $i\ne j$.
The needed cumulants for evaluating the profile cumulants in (\ref{profilecumulants}), and hence the $\alpha$-quantile modified profile score (\ref{modifiedprofilescore}), when $\psi=\beta_k$ or $\psi=\gamma_k$ are
factorized in terms involving separately the elements of the model matrices $X$ and $Z$, the link functions and cumulants of derivatives of $\log f(y_i;\mu_i,\phi_i)$. This particular structure of the cumulants is useful for computational implementation.

An {\tt R} implementation of the QBR confidence intervals for the regression models (\ref{linkmuphi}) is given in the function {\sffamily hoaci} available at the GitHub repository HOACI\footnote{\url{https://github.com/elianecpinheiro/HOACI}}, and two examples are presented, namely the Student-t and beta regression models, for which the cumulants are given in Examples 6 and 7 below.

\paragraph{Example 6.}{\it Symmetric and log-symmetric linear regression}.
Consider a heteroskedastic symmetric linear regression model as follows. Let $y_1,\ldots,y_n$ be independent random variables, each $y_i$ having a symmetric distribution ${\rm S}(\mu_i,\phi_i)$
with $\mu_i$ and $\phi_i$ as in (\ref{linkmuphi})
; see Example 5. 
The link functions are taken as $g(\mu_i)=\mu_i$ and $h(\phi_i)=\log(\phi_i)$. 
For this choice of link functions we have
$g^\prime(\mu)=1,$ $g^{\prime\prime}(\mu)=0,$
$h^\prime(\phi)=1/{\phi}$ and $h^{\prime\prime}(\mu)=-1/{\phi^2}.$
The needed cumulants are presented in Appendix \ref{ape:symregression}. 

A class of log-symmetric linear regression models for positive continuous responses is defined by assuming that $t_1,\ldots,t_n$ are such that 
$t_i=\exp(x_i^\top \beta)\xi_i^{\phi_i}$, where the $\xi_i$'s are independent and have a standard log-symmetric distribution with pdf $\xi^{-1}v(\xi^2)$, $\xi>0$
\citep{VanegasPaula2015, VanegasPaula2016}. An interesting feature of these models is that $\exp(x_i^\top \beta)$ is the median of $t_i$ and $\phi_i$ is
interpreted as a skewness parameter. 
Since $y_i = \log t_i \sim {\rm S}(x_i^\top \beta,\phi_i)$, the results above are also applicable 
for inference regarding the parameters of the log-symmetric linear regression models; see \citet[Sect. 4]{MedeirosFerrari2017}.

\paragraph{Example 7.}{\it Beta regression}. 
Beta regression models are widely applicable when the response variable is a continuous proportion \citep{FERRARICRIBARI, FERRARI}.
We consider the beta regression model that assumes that $y_i$ has a beta distribution with mean $\mu_i$ ($0 < \mu_i < 1$) and precision parameter $\phi_i$ ($\phi_i > 0$), and pdf
\begin{equation*}
\label{beta}
f(y_i; \mu_i, \phi_i) = \frac{\Gamma(\phi_i)}{\Gamma(\mu_i\phi_i) \Gamma((1-\mu_i)\phi_i)} y_i^{\mu_i\phi_i-1}(1-y_i)^{(1-\mu_i)\phi_i-1},\quad
0 < y_i < 1,
\end{equation*}
with $\mu_i$ and $\phi_i$ as in (\ref{linkmuphi}) 
\citep{SMITHSONVERKUILEN}.
The link functions are taken as the logit link for the mean, $g(\mu_i)=\log(\mu_i/(1-\mu_i))$, and the log link for the precision parameter, 
$h(\phi_i)=\log(\phi_i)$.
For this choice of link functions we have
$g^\prime(\mu)={1}/[\mu(1-\mu)],$
$g^{\prime\prime}(\mu)=(2\mu -1)/[\mu^2(1-\mu)^2],$
$h^\prime(\phi)={1}/{\phi}$ and
$h^{\prime\prime}(\mu)=-{1}/{\phi^2}.$
The needed cumulants are presented in Appendix \ref{ape:betaregression}.

\vspace{0.3cm}
It is worth mentioning that nonlinear regression models can be accomodated in the same framework. Let (\ref{linkmuphi}) be replaced by
$$
g(\mu_i ) = \eta_i = \eta(x_i, \beta) 
\quad {\rm and} \quad
h(\phi_i ) = \delta_i = \delta(z_i, \gamma ), 
$$
where $x_i$ and $z_i$ are known fixed vectors of dimensions $q^\prime$ and $m^\prime$ respectively, and 
$\eta(\cdot, \cdot)$ and $\delta(\cdot, \cdot)$ are allowed to be nonlinear functions in the second argument. 
Let $\cal X$ be the derivative matrix of $\eta = (\eta_1, \ldots , \eta_n)^\top$ with respect to $\beta^\top$. 
Analogously, let $\cal Z$ be the derivative matrix of $\delta = (\delta_1, \ldots , \delta_n)^\top$ with respect to $\gamma^\top$. 
In the linear case (\ref{linkmuphi}), ${\cal X} = X$ and ${\cal Z} = Z$. 
We assume that rank$({\cal X}) = q$ and rank$({\cal Z}) = m$ for all $\beta$ and $\gamma$. 
The needed cumulants of log-likelihood derivatives in the nonlinear case coincide with those for the linear case, with $X$ and $Z$ replaced by $\cal{X}$ and $\cal{Z}$, respectively.

\section{Monte Carlo simulation}\label{sec:MonteCarlo}

We present Monte Carlo simulation results to evaluate the finite-sample performance of confidence intervals based on the 
$\alpha$-quantile modified score (QBR confidence intervals). For comparison, we included results for Wald-type confidence intervals, i.e. confidence intervals that use the asymptotic normality of MLE (ML confidence intervals), and those that are constructed similarly to ML confidence intervals with median bias reduced estimates of all the parameters in place of MLEs (MBR confidence intervals). 
The simulations were implemented in {\tt R} \citep{R} and in the {\tt Ox} language \citep{DOORNIK}. 
The number of Monte Carlo replicates is 100,000. As initial guesses for the QBR confidence limits, we used the maximum likelihood estimates. The codes for the simulations are available at the GitHub repository HOACI.

The simulation results are shown in plots of `non-coverage discrepancy' of one-sided and two-sided confidence intervals, and the mean length of two-sided intervals. The non-coverage discrepancy of an interval is defined as the ratio of the non-coverage probability \ (evaluated via simulation) and 1 minus the \ nominal level. The \ non-coverage \ discrepancy is \ plotted against \ 1 minus the nominal level, i.e. the nominal non-coverage probability. \ Intervals with non-coverage discrepancy close to (greater than) 1 are those with coverage probability close to (smaller than) the nominal level. Since the intervals may not have the correct coverage probability, the mean length is plotted against the coverage probability and not the nominal level. \citet{PERETTI} suggest the use of mean length curves constructed this way to compare the `effectiveness' of different confidence regions.

We now list the scenarios for the simulation experiments.

\begin{figure}[!ht]
	\centering
		\includegraphics[width=4cm,height=4cm]{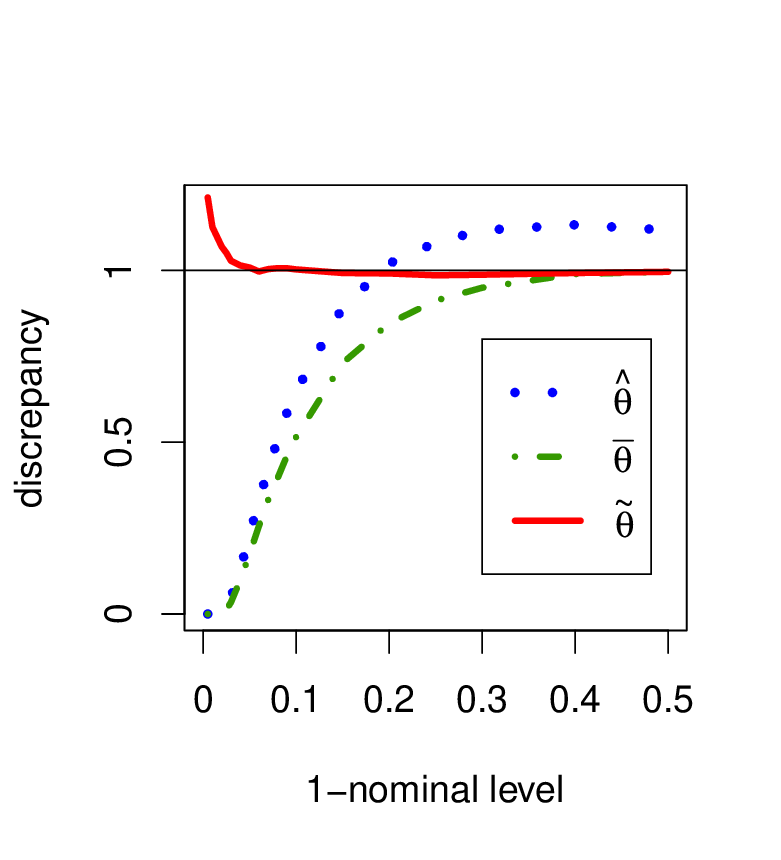}
		\includegraphics[width=4cm,height=4cm]{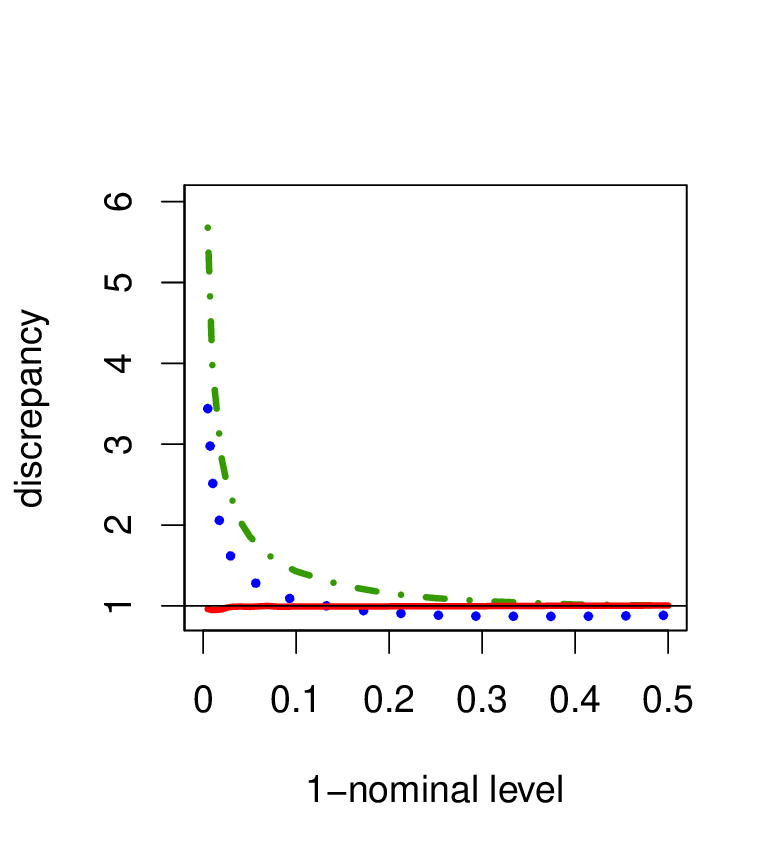}
		\includegraphics[width=4cm,height=4cm]{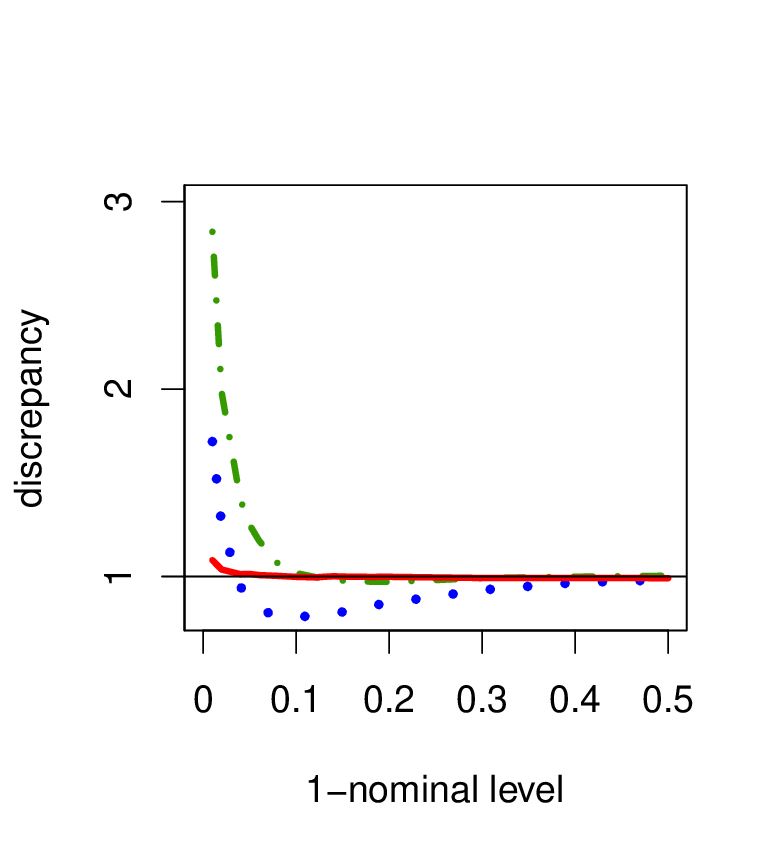}
		\includegraphics[width=4cm,height=4cm]{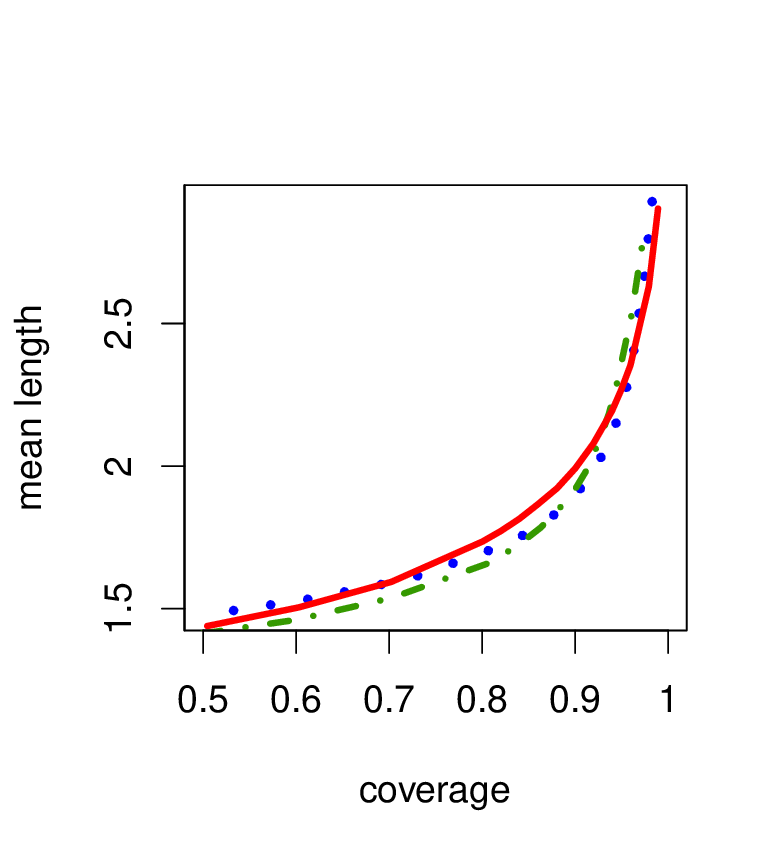}
	\caption{\footnotesize Plots of non-coverage discrepancy of the lower and upper one-sided confidence intervals (first and second plots) and non-coverage discrepancy  and mean length of two-sided confidence intervals (third and fourth plots); the legend indicates the parameter being estimated; hat, bar and tilde refer to ML, MBR, and QBR confidence intervals, respectively; exponential distribution. }
\label{fig:coveragediscrepancyexp}
\end{figure}

\paragraph{\bf Example 1 (cont.).} {\it One parameter exponential family.} 
We generated random samples of size $n=5$ from an exponential distribution with unit mean. Results are shown in Figure \ref{fig:coveragediscrepancyexp}.

\begin{figure}[!ht]
	\centering
		\includegraphics[width=4cm,height=4cm]{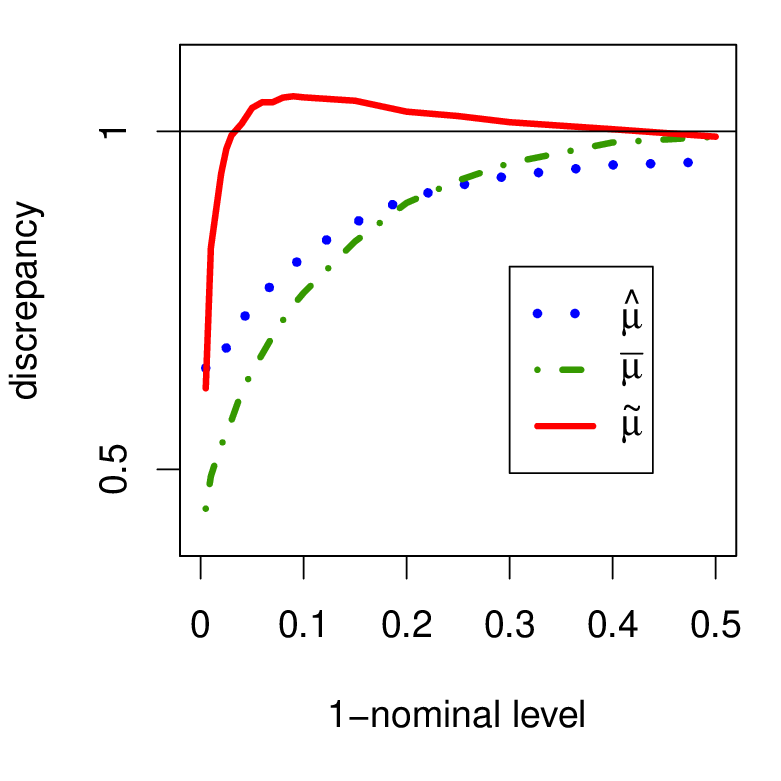}
		\includegraphics[width=4cm,height=4cm]{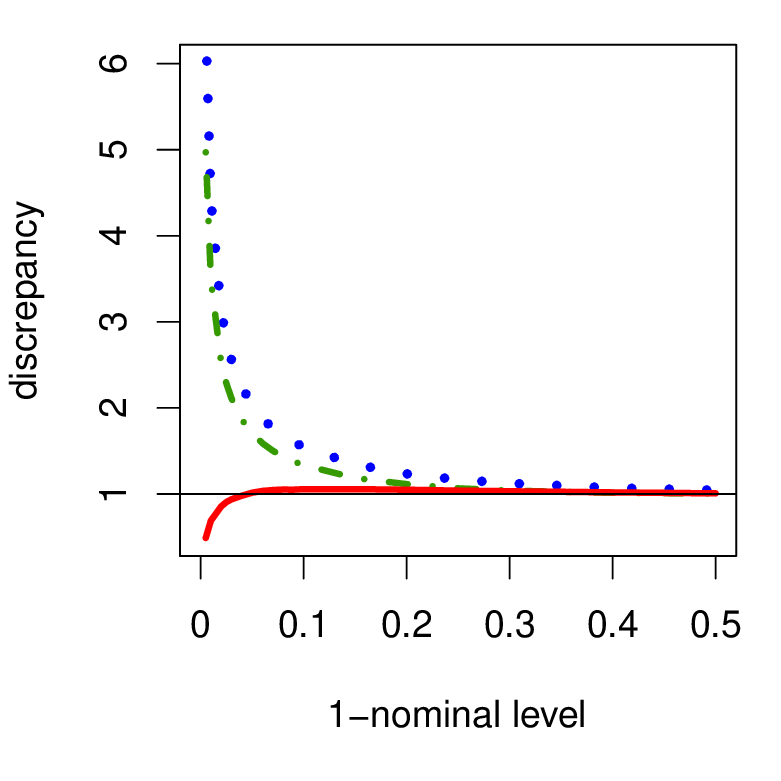}
		\includegraphics[width=4cm,height=4cm]{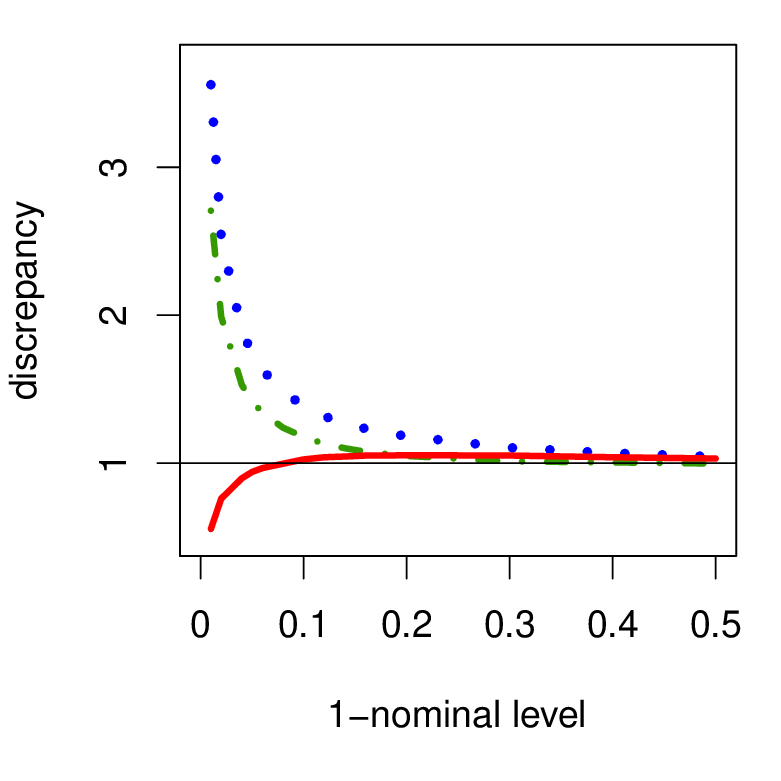}
 		\includegraphics[width=4cm,height=4cm]{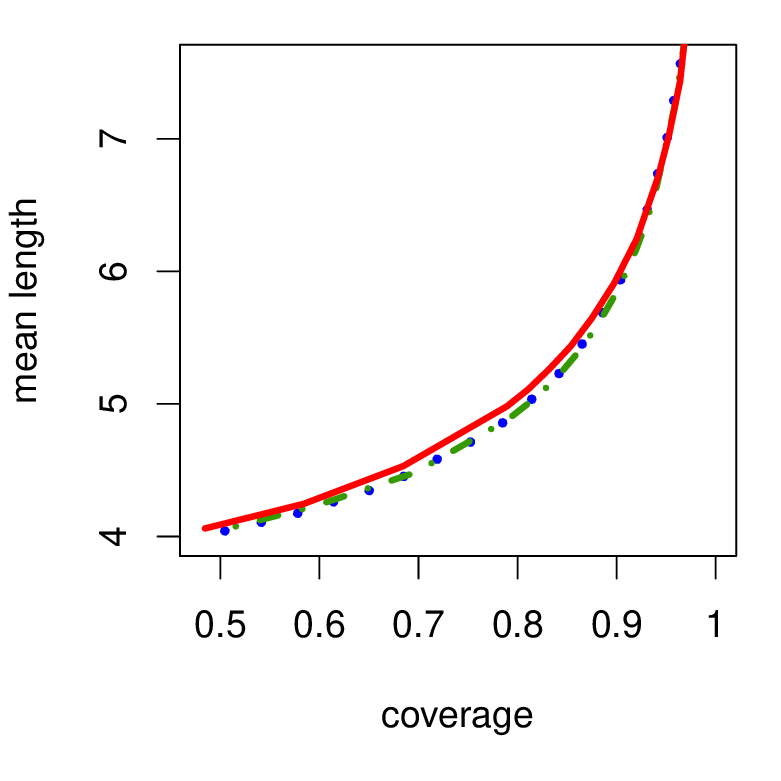}
\\
		\includegraphics[width=4cm,height=4cm]{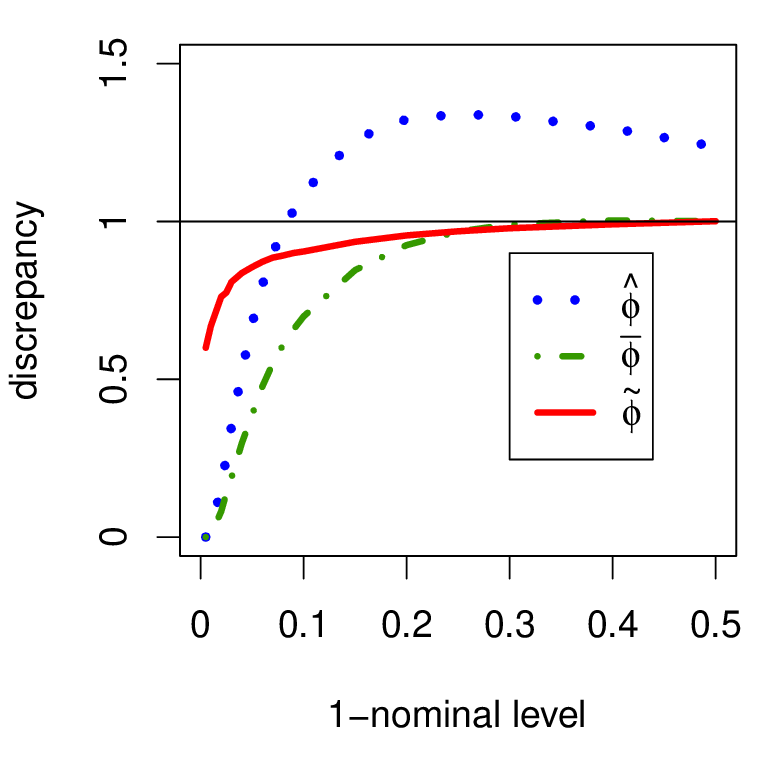}
		\includegraphics[width=4cm,height=4cm]{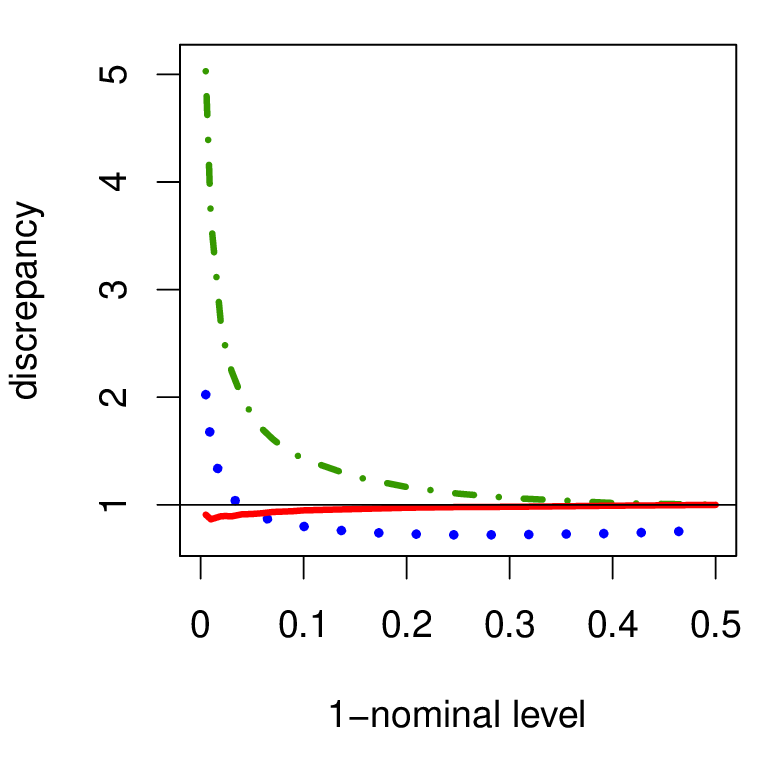}
		\includegraphics[width=4cm,height=4cm]{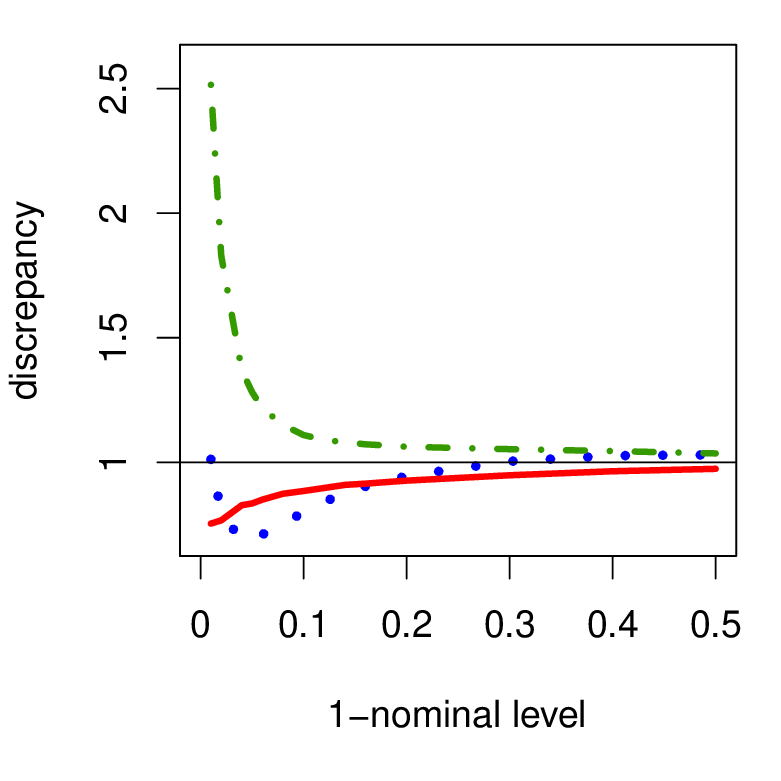}
 		\includegraphics[width=4cm,height=4cm]{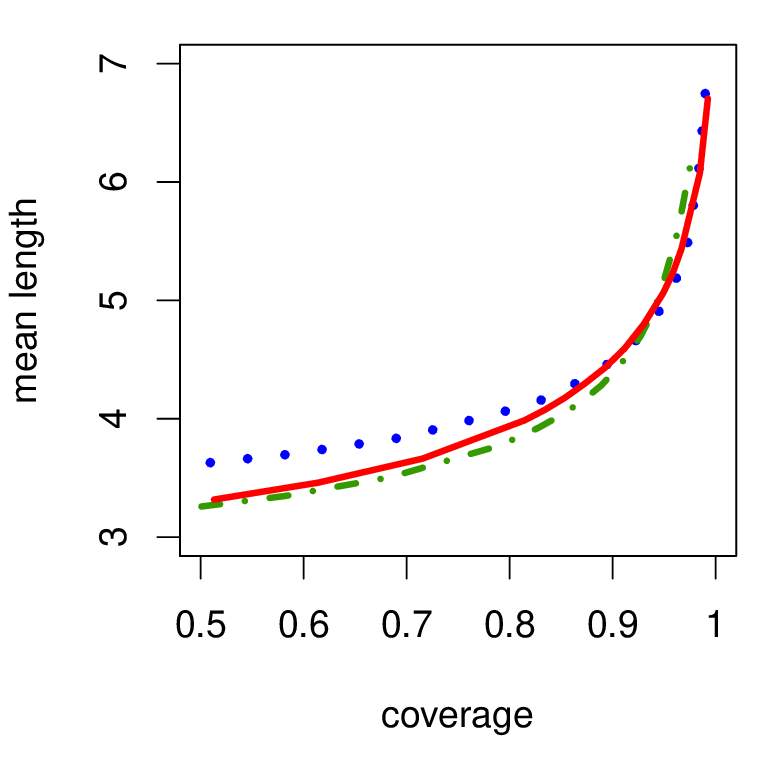}
\caption{\footnotesize Plots of non-coverage discrepancy of the lower and upper one-sided confidence intervals (first and second plots) and non-coverage discrepancy  and mean length of two-sided confidence intervals (third and fourth plots); the legend indicates the parameter being estimated; hat, bar and tilde refer to ML, MBR, and QBR confidence intervals, respectively; gamma distribution.} 
\label{fig:coveragediscrepancygamma}
\end{figure}

\vspace{0.1cm}
\paragraph{\bf Example 4 (cont.).}{\it Gamma distribution with mean $\mu$ and coefficient of variation $\phi^{-1/2}$, Gamma$(\mu,\phi)$.} 
We generated random samples of size $n=15$ from a  gamma distribution with $\mu=10$ and $\phi=3$. Results are given in Figure \ref{fig:coveragediscrepancygamma}.

\begin{figure}[!ht]
	\centering
		\includegraphics[width=4cm,height=4cm]{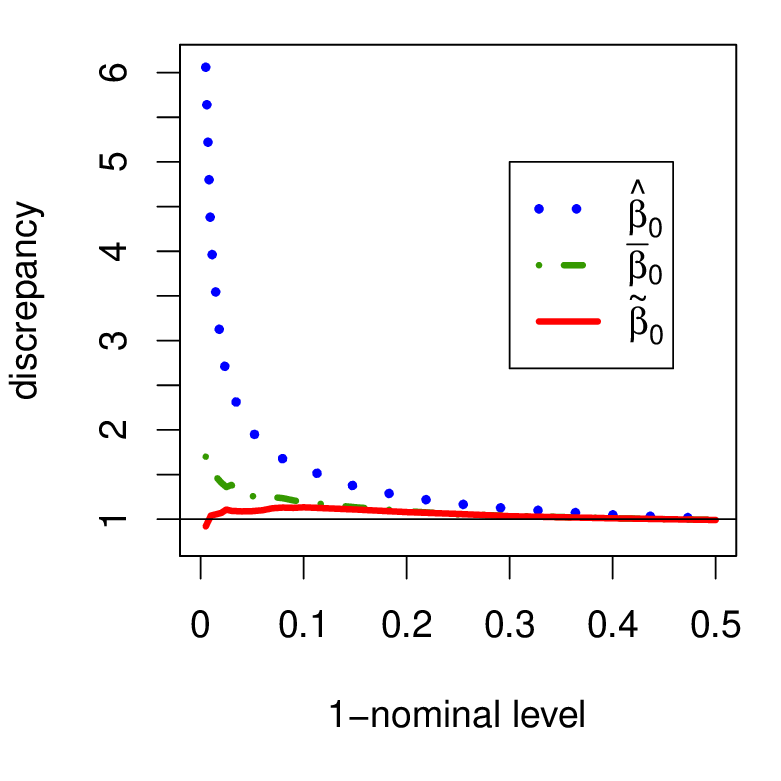}
		\includegraphics[width=4cm,height=4cm]{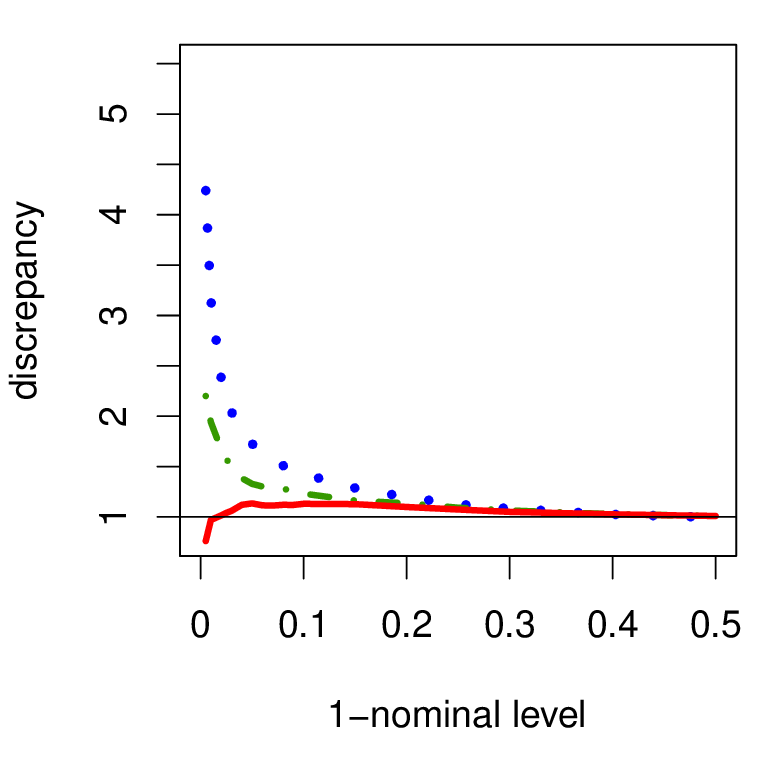}
		\includegraphics[width=4cm,height=4cm]{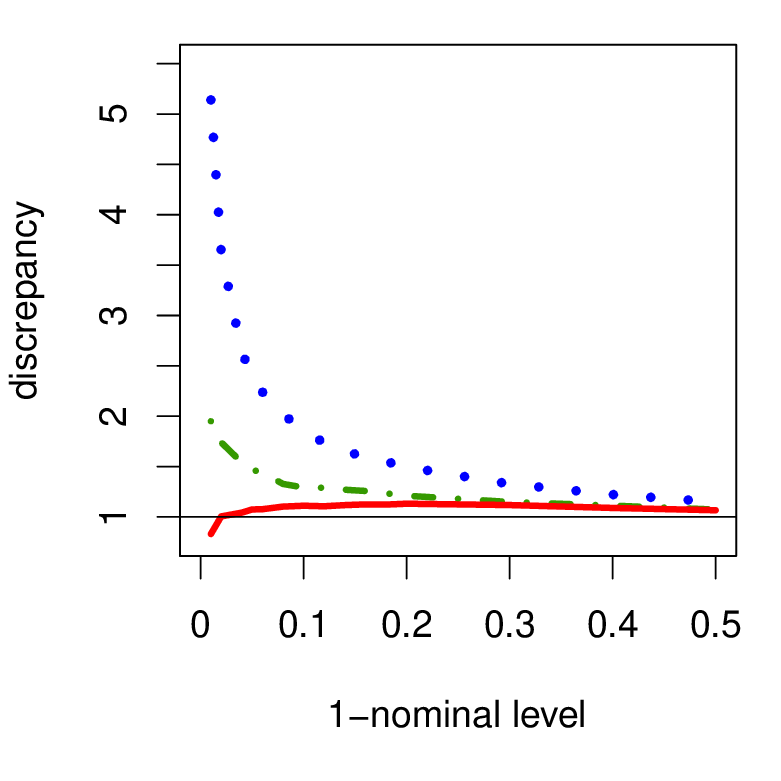}
		\includegraphics[width=4cm,height=4cm]{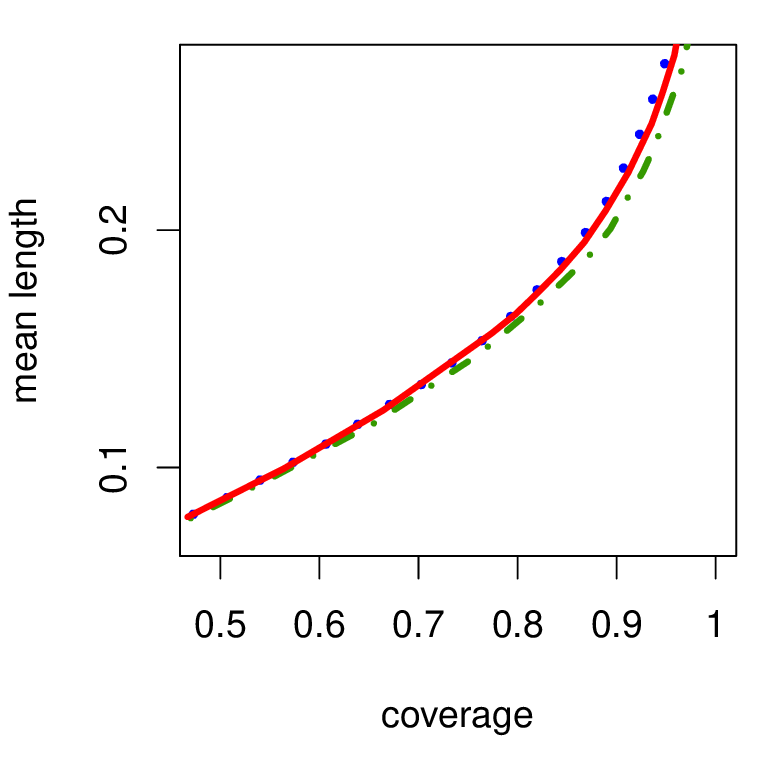}
		\\
		\includegraphics[width=4cm,height=4cm]{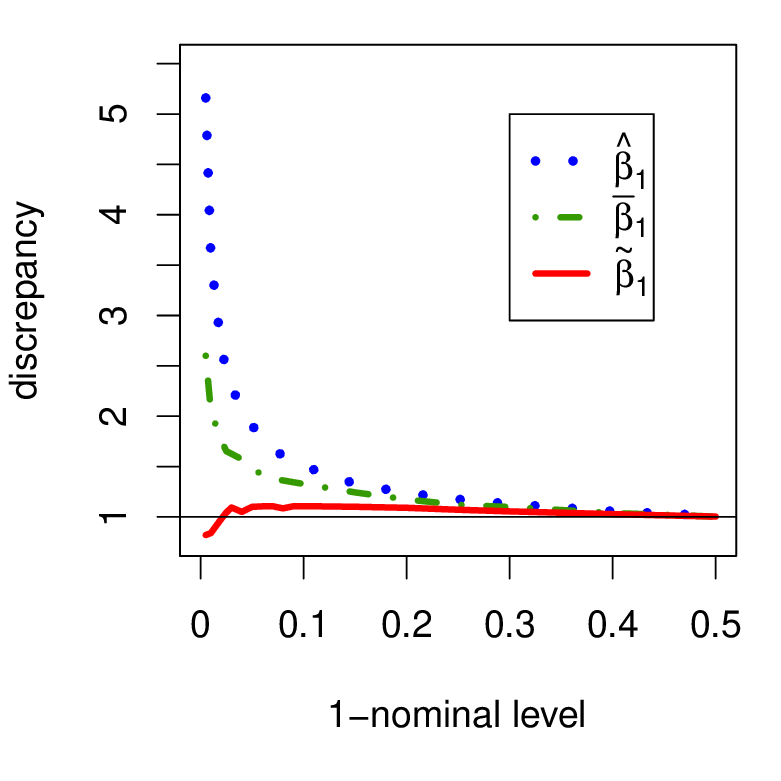}
		\includegraphics[width=4cm,height=4cm]{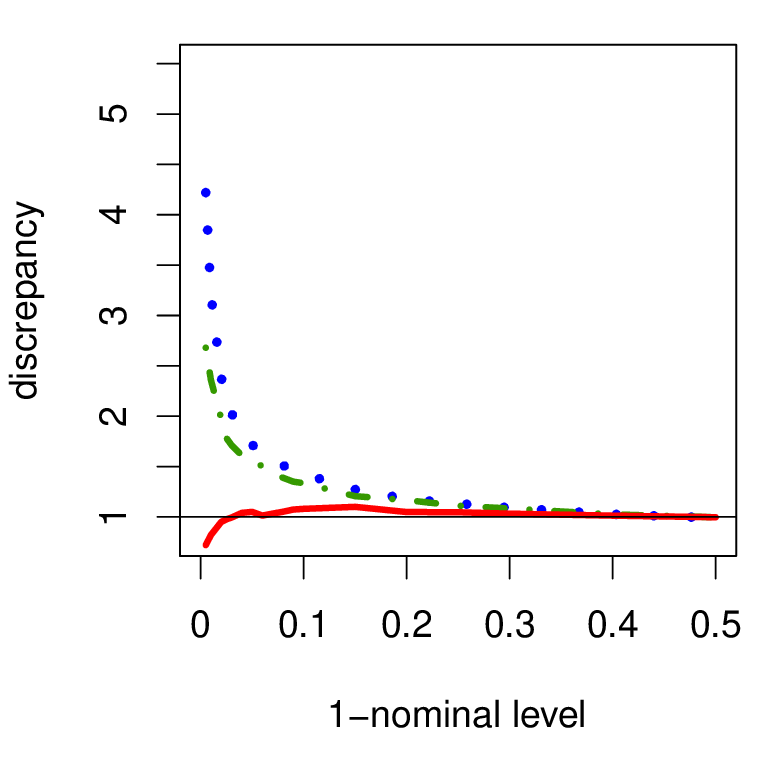}
		\includegraphics[width=4cm,height=4cm]{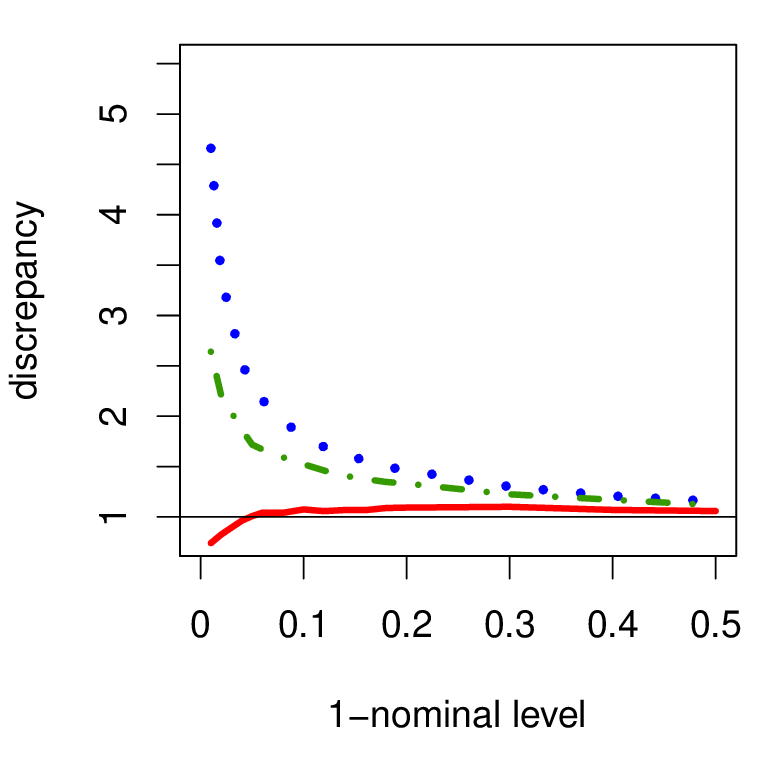}
		\includegraphics[width=4cm,height=4cm]{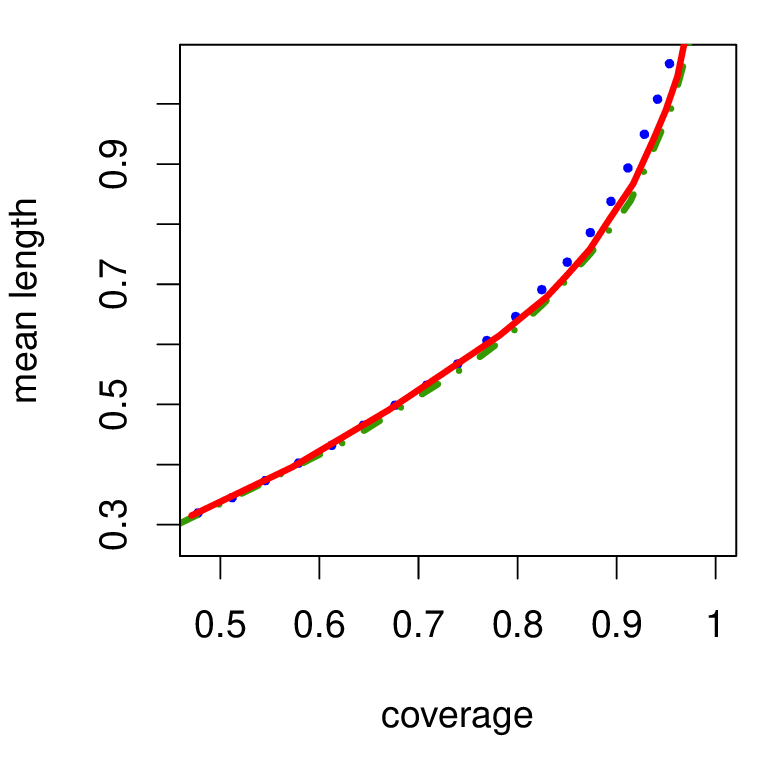}
		\\
		\includegraphics[width=4cm,height=4cm]{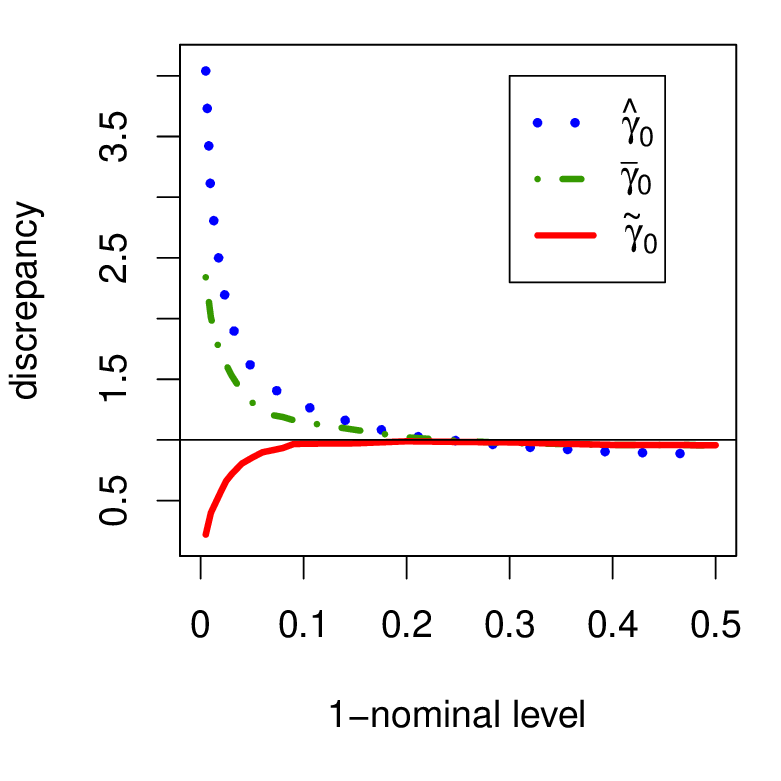}
		\includegraphics[width=4cm,height=4cm]{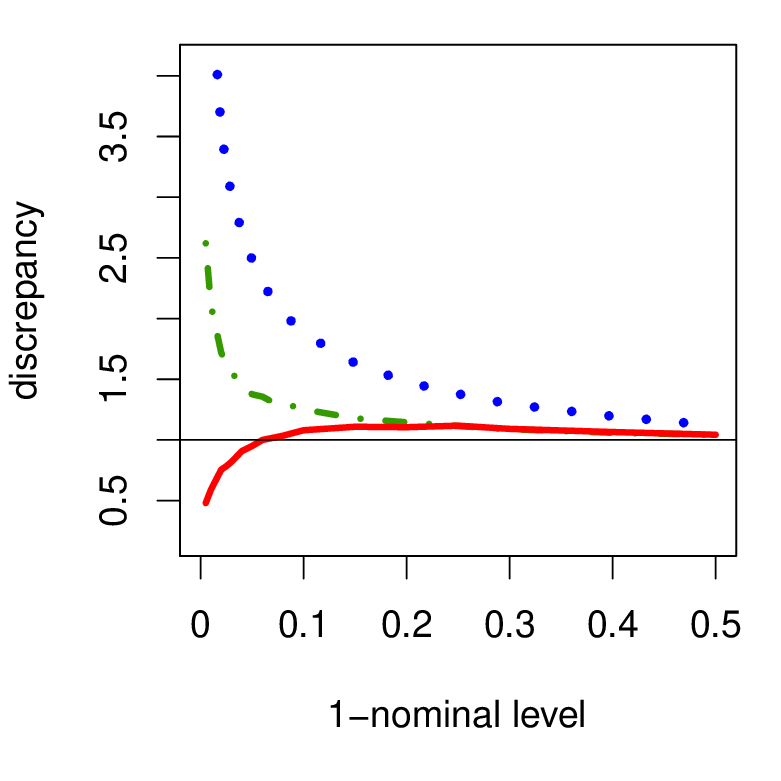}
		\includegraphics[width=4cm,height=4cm]{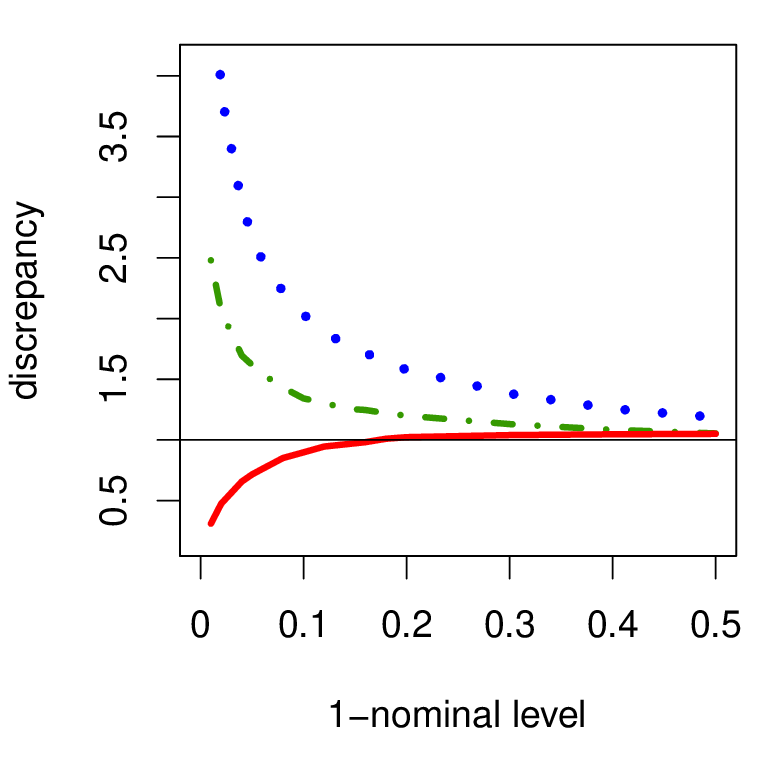}
		\includegraphics[width=4cm,height=4cm]{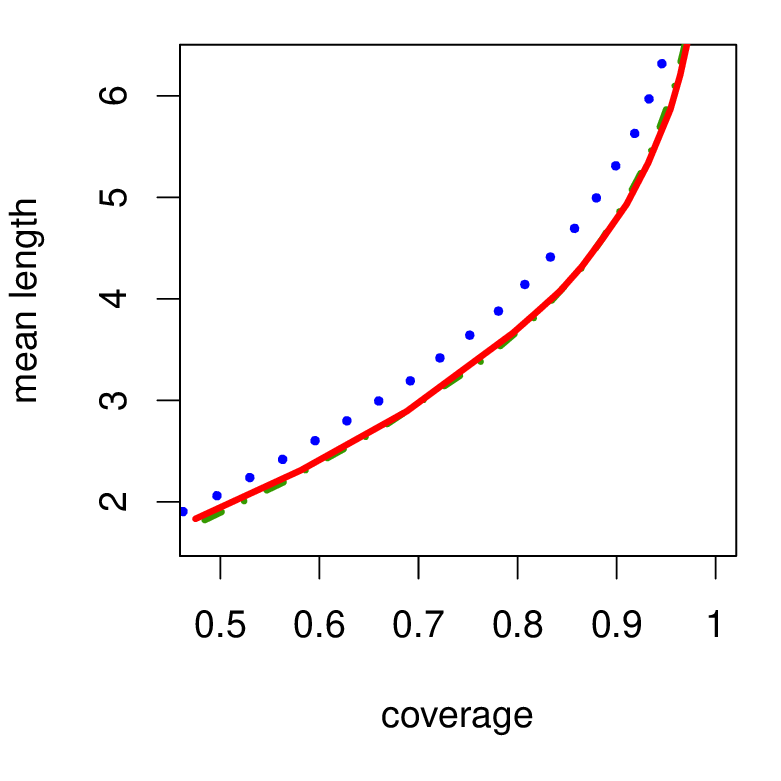}
		\\
		\includegraphics[width=4cm,height=4cm]{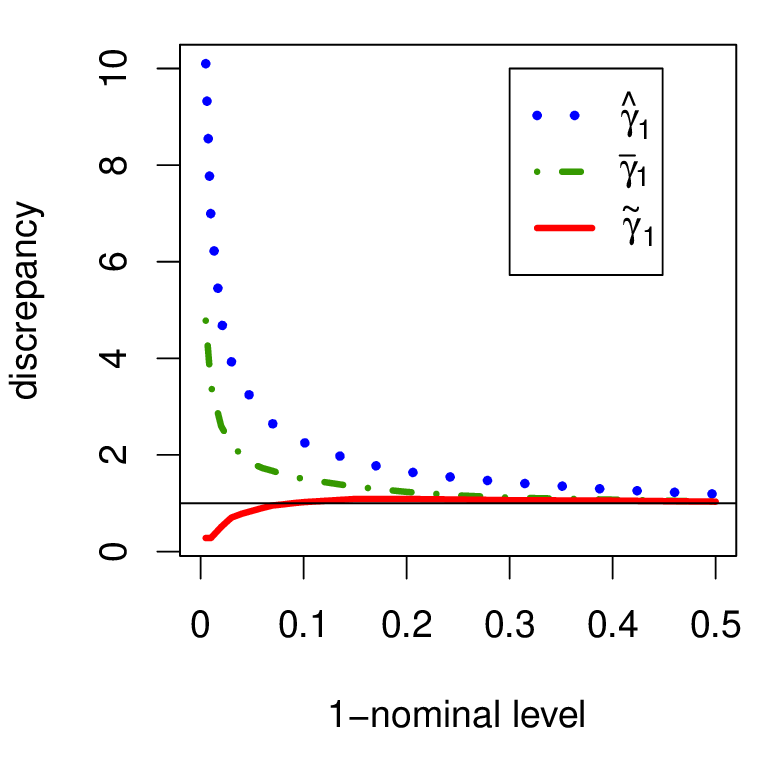}
		\includegraphics[width=4cm,height=4cm]{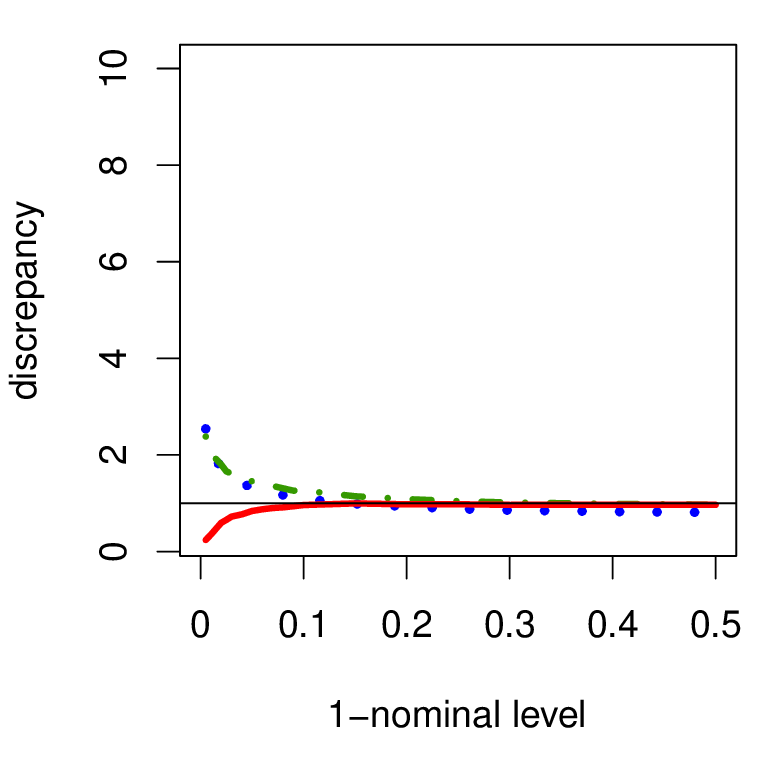}
		\includegraphics[width=4cm,height=4cm]{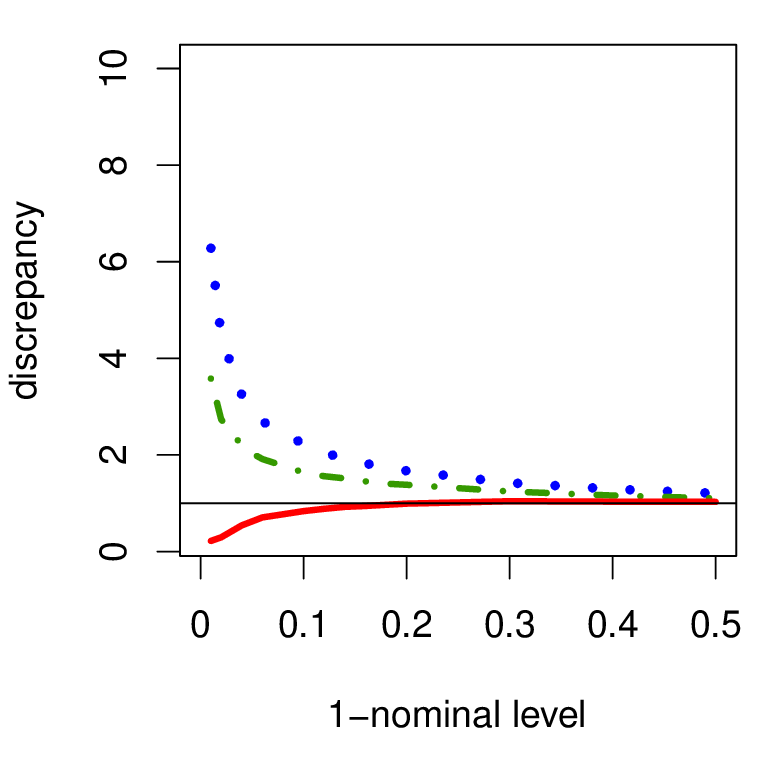}
		\includegraphics[width=4cm,height=4cm]{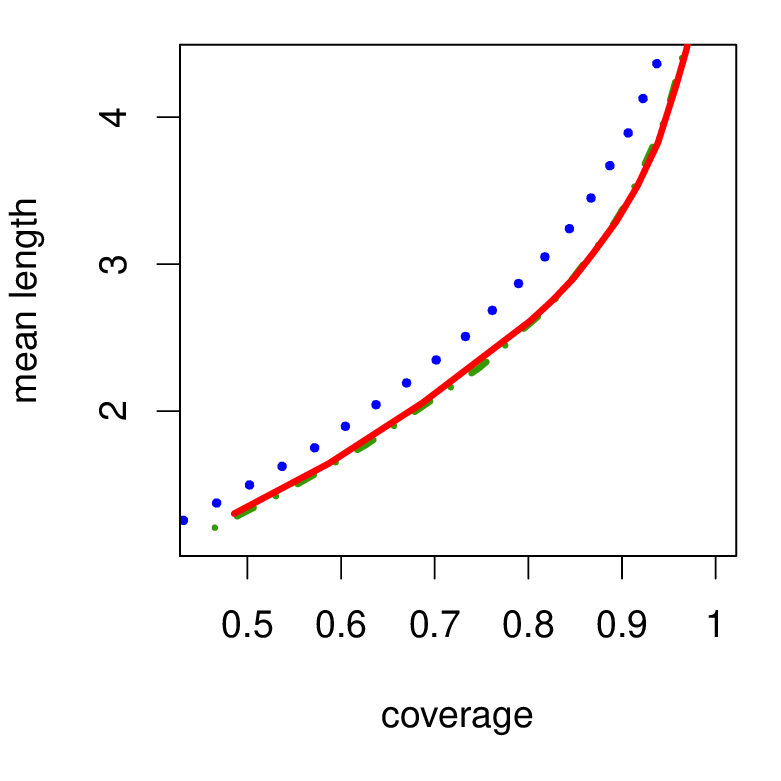}
	\caption{\footnotesize Plots of non-coverage discrepancy of the lower and upper one-sided confidence intervals (first and second plots) and non-coverage discrepancy  and mean length of two-sided confidence intervals (third and fourth plots); the legend indicates the parameter being estimated; hat, bar and tilde refer to ML, MBR, and QBR confidence intervals, respectively; beta regression.}
\label{fig:coveragediscrepancybetaregressionhetero}
\end{figure}

\vspace{0.1cm}
\paragraph{Example 7 (cont.).}{\it Beta regression}. 
We generated data from a beta regression model with mean and precision parameters defined, respectively, as
$\log[{\mu_i}/{(1-\mu_i)}] = \beta_{0} + \beta_{1}x_{i1}, 
\
\log(\phi_i) = \gamma_{0} + \gamma_{1}z_{i1}, \ i=1,\ldots,n,$
with $n=25$, $\beta_0=\beta_1=1$, $\gamma_0=1$, and $\gamma_1=2$. The values of the covariates $x$ and $z$ were drawn from uniform distributions in the intervals 
$(-1/2,1/2)$ and $(1,2)$, respectively. Results are seen in Figure \ref{fig:coveragediscrepancybetaregressionhetero}.

\vspace{1cm}
In Figures \ref{fig:coveragediscrepancyexp}-\ref{fig:coveragediscrepancybetaregressionhetero}, the blue, green, and red curves refer to ML, MBR, and QBR confidence intervals, respectively. From these figures the following conclusions may be drawn for the scenarios considered here. First, the non-coverage discrepancy tends to be much closer to 1 for the QBR confidence intervals, particularly for high nominal confidence levels, when compared with the ML confidence intervals. Second, in some cases the MBR confidence intervals partially correct the coverage probability of the ML confidence intervals but are clearly outperformed by the QBR confidence intervals proposed in this paper. Third, in some cases the lower and upper one-sided ML confidence intervals behave in opposite directions. For instance, for the gamma distribution (Figure \ref{fig:coveragediscrepancygamma}) the lower one-sided ML confidence interval for $\mu$ is conservative while the corresponding upper confidence interval has coverage probability smaller than the nominal level. This is expected because the asymptotic normality of MLEs does not account for the skewness of the MLEs in finite samples. This undesirable behavior does not occur when the QBR confidence intervals are employed. Finally, for each fixed coverage probability, the mean length of the QBR  two-sided confidence intervals tends to be similar to that of the ML two-sided confidence intervals. Overall, the simulations suggest that the method we propose performs considerably better than the usual approach, that employs the asymptotic normality of MLE when constructing confidence intervals.


\section{Applications}\label{sec:applications}

 We now present two applications using the {\sffamily hoaci} function available at the GitHub repository HOACI.

\subsection{Orange data}
\begin{figure}[!ht] 
	\centering
		\includegraphics[width=16cm,height=6cm]{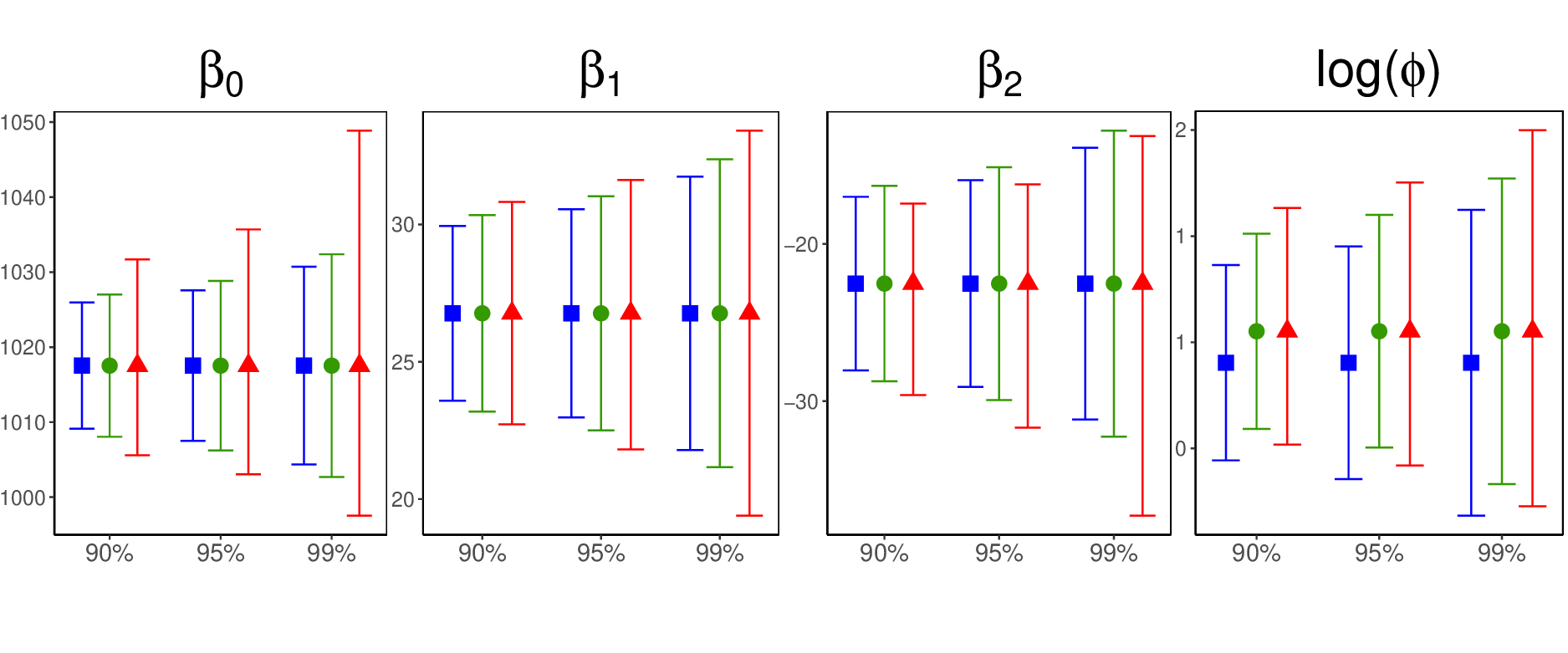}
	\caption{\footnotesize Plots of confidence intervals; blue square: ML, green circle: MBR, red triangle: QBR; orange data.}
\label{fig:CI_orange}
\end{figure}

The application considers the data set presented in Table 1 of \cite{MirhosseiniTan2010}. The data consist on 20 observations collected to investigate the effect of emulsion components on orange beverage emulsion properties. The response variable is the emulsion density, measured in miligrams per cubic centimeter, 
and the independent variables considered here are the amount of arabic gum ($x_1$) and of orange oil ($x_2$), both measured in grams per 10 grams. 
\citet[eq. (12)]{MedeirosFerrari2017} fitted a Student-t regression model with 3 degrees of freedom, 
\begin{eqnarray}
\label{symmodel}
\mu &=& \beta_0 + \beta_1 x_1 + \beta_2 x_2,
\end{eqnarray}
and constant dispersion parameter $\phi$. The point and interval estimates are given in Table \ref{tab:orange}. Figure (\ref{fig:CI_orange}) shows the confidence intervals for different confidence levels. The QBR confidence intervals are asymmetric, and they tend to be larger than the others, suggesting more uncertainty about the point estimates than the other confidence intervals. 

\begin{table}[!ht]
\begin{center}
\caption{\footnotesize Point estimates and confidence limits; 97.5\% confidence level for one-sided confidence interval, and 95\% confidence level for two-sided confidence interval; orange data} \label{tab:orange}
{\footnotesize
\begin{tabular}{lrrrrrrrrr}\hline
             &\multicolumn{2}{c}{Point estimates}& &\multicolumn{6}{c}{Confidence limits}\\
\cline{2-3}  \cline{5-10}
             & ML      & MBR                    & & \multicolumn{2}{c}{ML}& \multicolumn{2}{c}{MBR}& \multicolumn{2}{c}{QBR}   \\ 											
             &         &                        & & 2.5\%   &97.5\%   & 2.5\%   & 97.5\%  & 2.5\%   & 97.5\%   \\
\cline{2-3}  \cline{5-10}
\\[-1em]
$\beta_0$  &1017.5   &1017.5                    & &1007.5   &1027.6   &1006.2   &1028.8   &1003.0   &1035.7 \\
$\beta_1$  &26.8     &26.8                      & &23.0     &30.6     &22.6     &31.0     &21.8     &31.6 \\
$\beta_2$  &$-$22.5  &$-$22.5                   & &$-$29.1  &$-$15.9  &$-$29.9  &$-$15.1  &$-$31.7  &$-$16.2 \\
log(\ensuremath{\phi})&0.723    &0.841          & &0.285    &1.161    &0.403    &1.280    &0.335    &1.402 \\
\hline
\end{tabular}
} 
\end{center}       
\end{table}

\subsection{Reading skills}
The application considers the data set for assessing the contribution of non-verbal IQ to reading skills of dyslexic and non-dyslexic children \citep{SMITHSONVERKUILEN}. The data set comprises 44 observations and is included in the package {\tt betareg} \citep{BETAREG}. The independent variables are dyslexia ($x_1=$ $-1$ and $1$ for the control and dyslexic groups, respectively) and the nonverbal intelligent quotient ($x_2$, converted to $z$ scores), and the response variable is accuracy ($y$, the score on a test of reading accuracy). 
\cite{SMITHSONVERKUILEN} proposed the following beta regression model:
\begin{eqnarray}
\label{betamodel}
\log[{\mu}/{(1-\mu)}] &=& \beta_0 + \beta_1 \ x_1 + \beta_2 \ x_2 + \beta_3 \ x_1 x_2,
\nonumber \\
\log(\phi) &=& \gamma_0 + \gamma_1 \ x_1 + \gamma_2 \ x_2.
\end{eqnarray}
The point and interval estimates are given in Table \ref{tab:readingskills}. 

\begin{table}[!ht]
\begin{center}
\caption{\footnotesize Point estimates and confidence limits; 97.5\% confidence level for one-sided confidence interval, and 95\% confidence level for two-sided confidence interval; reading skills data} \label{tab:readingskills}
{\footnotesize
\begin{tabular}{lrrrrrrrrr}\hline
             &\multicolumn{2}{c}{Point estimates}& &\multicolumn{6}{c}{Confidence limits}\\
\cline{2-3}  \cline{5-10}
             & ML      & MBR                     & & \multicolumn{2}{c}{ML}& \multicolumn{2}{c}{MBR}& \multicolumn{2}{c}{QBR}   \\
             &         &                         & & 2.5\%      &97.5\%    & 2.5\%      & 97.5\%    & 2.5\%   & 97.5\%   \\
						
\cline{2-3}  \cline{5-10}
\\[-1em] 											
$\beta_0$  &1.12    &1.11                    & &0.84    &1.40    &0.82    &1.40    &0.73    &1.42    \\
$\beta_1$  &$-$0.74 &$-$0.73                 & &$-$1.02 &$-$0.46 &$-$1.02 &$-$0.44 &$-$1.04 &$-$0.35 \\
$\beta_2$  &0.49    &0.47                    & & 0.23   &0.75    &0.19    &0.75    &0.00    &0.86    \\
$\beta_3$  &$-$0.58 &$-$0.57                 & &$-$0.84 &$-$0.32 &$-$0.84 &$-$0.29 &$-$0.95 &$-$0.04 \\
$\gamma_0$ &3.30     &3.11                      & &2.87  &3.74  &2.67  &3.55  &2.41    &3.60 \\
$\gamma_1$ &1.75     &1.69                      & &1.23  &2.26  &1.18  &2.21  &0.92    &2.27 \\
$\gamma_2$ &1.23     &1.06                      & &0.71  &1.75  &0.53  &1.60  &$-$0.33 &2.38 \\
\hline
\end{tabular}
} 
\end{center}
\end{table}

\begin{figure}[!ht]
	\centering
		\includegraphics[width=16cm,height=6cm]{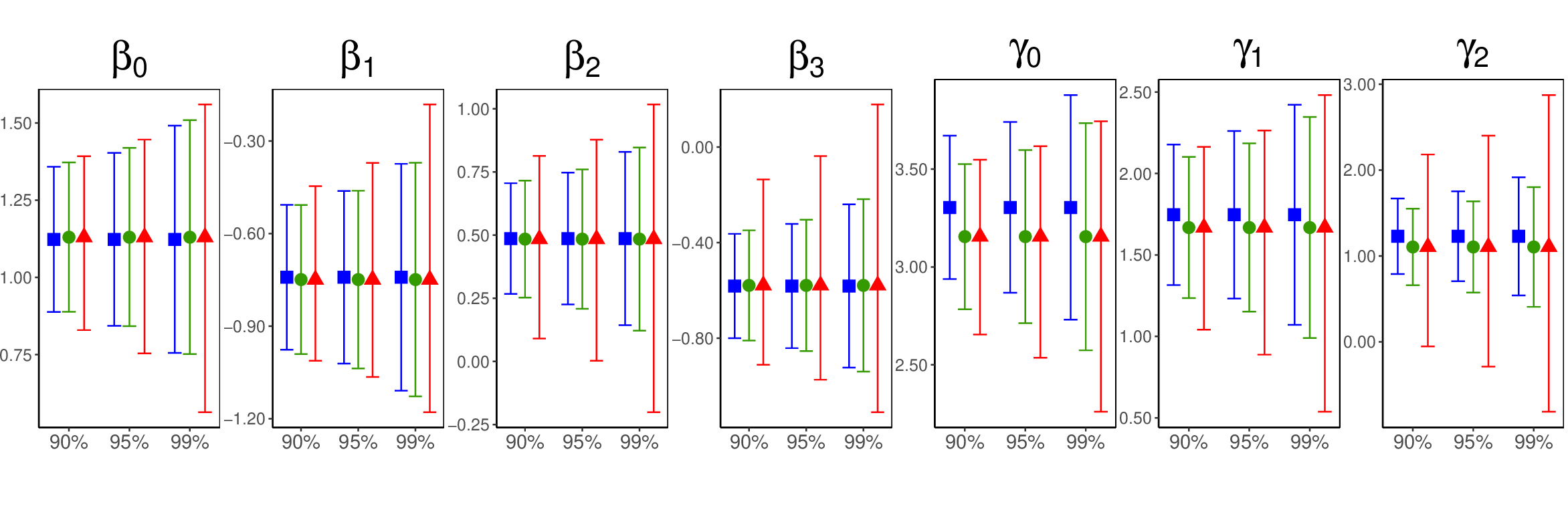}
	\caption{\footnotesize Plots of confidence intervals; blue square: ML, green circle: MBR, red triangle: QBR; reading skills data.}
\label{fig:CI_reading}
\end{figure}

The QBR confidence intervals are quite different from the others in some cases. For instance, the 95\%  
confidence intervals for $\gamma_2$ are $(0.71, 1.75)$ (ML), $(0.53, 1.60)$ (MBR), and $(-0.33, 2.38)$ (QBR). 
Hence, the QBR confidence interval does not provide evidence of effect of non-verbal IQ in the variability of reading accuracy unlike the other confidence intervals. In other words, the QBR confidence intervals change the inferential conclusions. Figure (\ref{fig:CI_reading}) shows the confidence intervals for different confidence levels. Considerable differences in the length and shape (asymmetry) of the intervals are observed, suggesting that the normal approximation for the MLEs is not accurate. 

To further investigate the accuracy of the different confidence intervals in the present scenario, we conduct a simulation experiment using the beta regression model (\ref{betamodel}) with $n=44$ and the values of the parameters taken as the MLEs computed with the reading skills data set. The simulation experiment includes the confidence intervals obtained from the signed log likelihood ratio statistic ($r$) and its third-order modification, $r^*$ \citep{BARNDORFF1986}; for a recent review on inference based on $r^*$ and an associate {\tt R} package, see \citet{PIERCEBELLIO}. Here, we considered Skovgaard's approximation \citep{SKOVGAARD2001} to $r^*$ as implemented in \citet{PIERCEBELLIO}. The coverage of the confidence intervals at nominal levels $90\%$, $95\%$, and $99\%$ are presented in Table \ref{tab:readingskillssimulation}. It is clear from the table that the ML confidence intervals are anti-conservative (undercover) and the QBR confidence intervals proposed here are slightly conservative but cover the true values of the parameters with probability much closer to the nominal levels than the ML confidence intervals. The confidence intervals for $\gamma_0$ based on $r$ and on the MLE are considerably under-conservative. Additionally, the performances of the QBR confidence intervals and those based on $r^*$ are comparable, with slight superiority of the  latter in this particular simulation scenario. The slight superiority of the $r^*$ confidence intervals might be due to the fact that they control the relative error, rather than the absolute error, as the proposed method. Additional Monte Carlo investigation in different scenarios comparing these different approaches would be helpful but it is beyond the scope of this article.

\begin{table}[!ht]
\begin{center}
{\footnotesize
\caption{\footnotesize Coverage of confidence intervals (in percentage); reading skills data scenario} 
\label{tab:readingskillssimulation}
\begin{tabular}{r r r r r r r r r}
  \hline
  & & $\beta_0$ & $\beta_1$ & $\beta_2$ & $\beta_3$ & $\gamma_0$ & $\gamma_1$ & $\gamma_2$ \\ 
  \hline
   90\% & ML   & 86.9 & 87.0 & 85.0 & 84.8 & 75.8 & 82.2 & 78.7 \\ 
        & MBR  & 88.1 & 88.1 & 87.3 & 87.1 & 86.4 & 83.2 & 82.6 \\
        & QBR  & 90.6 & 90.9 & 90.5 & 90.3 & 93.0 & 90.4 & 92.6 \\
				& $r$  & 88.7 & 88.7 & 87.6 & 87.3 & 78.8 & 86.9 & 86.5 \\
        & $r^*$& 89.3 & 89.2 & 89.6 & 89.3 & 89.9 & 88.9 & 89.1 \\
				\hline
   95\% & ML   & 92.5 & 92.5 & 91.2 & 91.0 & 83.5 & 88.9 & 86.6 \\ 
        & MBR  & 93.4 & 93.5 & 93.0 & 92.9 & 92.3 & 89.9 & 88.7 \\
        & QBR  & 96.2 & 96.1 & 95.8 & 95.6 & 96.9 & 96.1 & 97.1 \\
        & $r$  & 94.1 & 94.0 & 93.4 & 93.3 & 86.4 & 92.7 & 92.4 \\
        & $r^*$& 94.7 & 94.6 & 94.7 & 94.8 & 94.8 & 94.5 & 94.3 \\				
				\hline
   99\% & ML   & 97.7 & 97.5 & 97.3 & 97.3 & 93.1 & 96.2 & 94.4 \\ 
        & MBR  & 98.2 & 98.0 & 98.1 & 98.0 & 97.9 & 96.6 & 96.2 \\ 
        & QBR  & 99.4 & 99.4 & 99.3 & 99.4 & 99.4 & 99.7 & 99.7 \\
				& $r$  & 98.8 & 98.6 & 98.3 & 98.4 & 95.8 & 98.3 & 98.1 \\
        & $r^*$& 98.9 & 98.8 & 98.9 & 98.9 & 98.9 & 98.9 & 98.7 \\
   \hline
\end{tabular}
}
\end{center}
\end{table}

\section{Concluding remarks}\label{sec:conclusion}

We derived highly accurate confidence intervals for a scalar parameter of interest possibly in the presence of a vector of nuisance parameters in general parametric families. The proposed confidence limits are computed from modified score equations and possess desirable properties: they are equivariant under interest-respecting reparameterizations, they account for skewness and kurtosis of the score function, and they 
do not require computer-intensive resampling methods.

The $\alpha$-quantile modified (profile) score functions in (\ref{modifiedscore}) and (\ref{modifiedprofilescore}) generalize the median-modified (profile) score functions derived in \citet[eq. (1) and (8)]{PAGUI}. Hence, by setting $\alpha=0.5$, they provide third order median unbiased estimates for a scalar parameter in the presence or absence of nuisance parameters. \citet{PAGUI} also present an alternative method to obtain third order median unbiased estimates in the multiparameter setting by a suitable modification in the score vector (not the profile score); see eq. (10) in that paper. We have not explored this approach here, and this topic in the context of modified quantile score functions is an open route for subsequent research.

A relevant topic for future work is a comprehensive numerical comparison of the proposed confidence limits with alternatives, such as confidence limits based on the signed log likelihood ratio statistic, $r$, its third-order modification, $r^*$, the score statistic, the gradient statistic, the Bartlett-corrected likelihood ratio statistic \citep{CORDEIROCRIBARI}, the Bartlett-type corrected score statistic \citep{CORDEIROFERRARI},  the Bartlett-type corrected gradient statistic \citep{VARGAS}, and
modified profile likelihoods.

\paragraph{Acknowledgements.} We thank two anonymous reviewers for fruitful comments and suggestions that improved the paper. We also thank Euloge Clovis Kenne Pagui for sharing his R code for evaluating the median bias-corrected estimate of the skewness parameter of the skew-normal distribution. The second author gratefully acknowledges funding provided by Conselho Nacional de Desenvolvimento Científico e Tecnológico--CNPq (Grant No. 305963-2018-0).

\small
\begin{appendices}
\section{Proof of (\ref{modifiedscorefda})}
\label{ape:approximationorder}

The Edgeworth expansion for the cumulative distribution function (cdf) of a standardized sum of independent random variables, say $S_n^*$,  is
$$
F_{S_n^*}(x)=\Phi(x)-\phi(x)h(x) + O(n^{-3/2})
$$
where $\Phi(\cdot)$ and $\phi(\cdot)$ are the cdf and the pdf of the standard normal distribution and
\begin{eqnarray*}
h(x)&=& 
\frac{1}{6}\frac{\kappa_3}{\kappa_2^{3/2}}(x^2-1) 
+ \frac{1}{24}\frac{\kappa_4}{\kappa_2^{2}}(x^3-3x)
+\frac{1}{72}\frac{\kappa_3^2}{\kappa_2^{3}}(x^5-10x^3+15x);
\end{eqnarray*}
\citet[eq. (10.14)]{PACE}.
Let 
$$
h_1(x)=\frac{1}{6}\frac{\kappa_3}{\kappa_2^{3/2}}(x^2-1), \ \ 
h_2(x)=\frac{1}{24}\frac{\kappa_4}{\kappa_2^{2}}(x^3-3x)
-\frac{1}{72}\frac{\kappa_3^2}{\kappa_2^{3}}(4x^3-10x), \ \ 
h_3(x)=\frac{1}{72}\frac{\kappa_3^2}{\kappa_2^{3}}(x^5-6x^3+5x).
$$
Note that $h(x)= h_1(x) + h_2(x) + h_3(x)$; $h_1(x)$ is of order $O(n^{-1/2})$, and $h_2(x)$ and $h_3(x)$ are of order $O(n^{-1})$. 

From (\ref{modifiedscore}) and applying the Edgeworth expansion given above to the cdf of 
$U(\theta)/\sqrt{\kappa_2}$ we have
\begin{eqnarray*}
P_\theta\left(\widetilde U(\theta)\leq 0 \right) 
&=& 
P_\theta
\left(\frac{U(\theta)}{\sqrt{\kappa_2}}  \leq  u_\alpha + h_1(u_\alpha) + h_2(u_\alpha)\right)
\\  &=&
\Phi(u_\alpha + h_1(u_\alpha) + h_2(u_\alpha)) - \phi(u_\alpha + h_1(u_\alpha) + h_2(u_\alpha))
\ h(u_\alpha + h_1(u_\alpha) + h_2(u_\alpha)) \\
&&+ O(n^{-3/2}).
\end{eqnarray*}
Now, applying a Taylor series expansion and using the fact that $\phi^\prime(u_\alpha)=-u_\alpha\phi(u_\alpha)$, where prime denotes derivative with respect to $u_\alpha$, we have after some algebra
\begin{eqnarray*}
P_\theta\left(\widetilde U(\theta)\leq 0 \right) 
&=& 
\Phi(u_\alpha) - \phi(u_\alpha) \ h(u_\alpha)
+ \left(\Phi(u_\alpha) - \phi(u_\alpha) \ h(u_\alpha)\right)^\prime \ (h_1(u_\alpha) + h_2(u_\alpha))
\\&&
+ \frac{1}{2}\left(\Phi(u_\alpha) - \phi(u_\alpha) \ h(u_\alpha)\right)^{\prime\prime} \ (h_1(u_\alpha) + h_2(u_\alpha))^2
+O(n^{-3/2})
\\  &=&
\alpha 
- \phi(u_\alpha) \ h_3(u_\alpha)
- \phi^\prime(u_\alpha) \ h_1^2(u_\alpha) 
- \phi(u_\alpha) \ h_1^\prime(u_\alpha) \ h_1(u_\alpha)
+ \frac{1}{2} \phi^\prime(u_\alpha) \ h_1^2(u_\alpha)\\
&&+O(n^{-3/2})
\\  &=&
\alpha 
+ \phi(u_\alpha) 
\left(
-h_3(u_\alpha)
+ \frac{u_\alpha}{2} \ h_1^2(u_\alpha) 
- h_1^\prime(u_\alpha) \ h_1(u_\alpha)
\right)
+O(n^{-3/2})
\\  &=&
\alpha 
+ \phi(u_\alpha) \frac{1}{72}\frac{\kappa_3^2}{\kappa_2^{3}} 
\left(
-(u_\alpha^5 - 6 u_\alpha^3 + 5u_\alpha)
+ u_\alpha (u_\alpha^2-1)^2 
- 4 u_\alpha (u_\alpha^2-1)
\right)
+O(n^{-3/2})
\\  &=&
\alpha + O(n^{-3/2}).
\end{eqnarray*}
\section{Symmetric regression: cumulants}
\label{ape:symregression}

The first and second order derivatives of $\log f(y; \mu, \phi)$ are given by 
$U_{\mu}=-\phi^{-1}s^{(1)},$
$U_{\phi}= -\phi^{-1}(1 + s^{(1)}\epsilon),$
$U_{\mu\mu}= \phi^{-2}s^{(2)},$
$U_{\phi \phi}= \phi^{-2} (1 + 2s^{(1)}\epsilon + s^{(2)}\epsilon^2),$
and
$U_{\mu\phi}=\phi^{-2}  (s^{(1)} + s^{(2)}\epsilon),$
where $s^{(r)}=\dd^r s(\epsilon)/\dd \epsilon^r$ with $s(\epsilon)=\log v(\epsilon^2)$ and $\epsilon=(y-\mu)/\phi$.
From these derivatives, we obtain the following  cumulants:

\begin{eqnarray*}
\kappa_{\mu,\mu}             &=& \delta_{20000}/\phi^{2},
\ \
\kappa_{\phi,\phi}            =  (\delta_{20002}-1)/\phi^{2},
\ \
\kappa_{\mu,\mu,\phi}         =  2 \delta_{11001}/\phi^{3},
\ \
\kappa_{\phi,\phi,\phi}       =  2 (\delta_{11003}+1)/\phi^{3},
\\
\kappa_{\mu,\mu,\mu,\mu}     &=& (\delta_{40000}-3\delta_{20000}^2)/\phi^{4},
\ \
\kappa_{\mu,\mu,\mu,\phi}     =  (\delta_{30000}+\delta_{40001})/\phi^{4},
\\
\kappa_{\mu,\mu,\phi,\phi}   &=&  (2 \delta_{30001}+\delta_{40002}-\delta_{01000} \delta_{01002})/\phi^{4},
\ \
\kappa_{\mu,\phi,\phi,\phi}   =  (\delta_{40001}+3\delta_{30001}+3\delta_{20001})/\phi^{4},
\\
\kappa_{\phi,\phi,\phi,\phi} &=& (\delta_{40004} + 4 \delta_{30003}+12 \delta_{20002} - 3 \delta_{20002}^2-6)/\phi^{4},
\ \
\kappa_{\mu,\mu\phi}          =  -(\delta_{11001} - \delta_{01000})/\phi^{3},
\\
\kappa_{\phi,\phi\phi}       &=& (4 \delta_{01002} + \delta_{00103} - 2)/\phi^{3},
\ \
\kappa_{\phi,\mu\mu}          =  \delta_{00101}/\phi^{3},
\\
\kappa_{\mu,\phi}            &=& \kappa_{\mu,\mu,\mu} = \kappa_{\mu,\mu\mu} = \kappa_{\mu,\phi,\phi} = \kappa_{\mu,\phi\phi} = \kappa_{\phi,\mu\phi} =0.
\end{eqnarray*}
\section{Beta regression: cumulants}
\label{ape:betaregression}

The first and second order derivatives of $\log f(y; \mu, \phi)$ are given by
$$
U_{\mu}=\phi \left( y^* - \mu^* \right),
\quad\quad
U_{\phi}=\mu(y^* - \mu^*) + (y^\dag - \mu^\dag),
$$
$$
U_{\mu\mu}=-\phi^2\Bigl( \Psi^{(1)}(\mu \phi ) + \Psi^{(1)}((1-\mu )\phi ) \Bigr), 
\quad\quad
U_{\phi \phi}= -\mu^2 \Psi^{(1)}(\mu \phi) + \Psi^{(1)}(\phi), 
$$
$$
U_{\mu\phi}= y^*-\mu^*- \phi (\mu \Psi^{(1)}(\mu \phi) - (1-\mu)\Psi^{(1)}((1-\mu)\phi) ),
$$
where
$y^* = \log(y/(1-y))$, 
$y^{\dag}= \log(1-y)$, 
$\mu^*     = \Psi^{(0)}(\mu\phi) - \Psi^{(0)}((1-\mu)\phi)$, and
$\mu^\dag  = \Psi^{(0)}((1-\mu)\phi)-\Psi^{(0)}(\phi)$. 
From these derivatives, we obtain the following cumulants:
\begin{eqnarray*}
\kappa_{\mu,\mu}             &=&  \phi^2 (      \Psi^{(1)}(\mu\phi) +            \Psi^{(1)}((1-\mu)\phi) ),
\\
\kappa_{\mu,\phi}            &=&  \phi   (\mu   \Psi^{(1)}(\mu\phi) - (1-\mu)  \Psi^{(1)}((1-\mu)\phi) ),
\\
\kappa_{\phi,\phi}           &=&          \mu^2 \Psi^{(1)}(\mu\phi) + (1-\mu)^2\Psi^{(1)}((1-\mu)\phi) - \Psi^{(1)}(\phi),
\\
\kappa_{\mu,\mu,\mu}         &=&  \phi^3 (      \Psi^{(2)}(\mu\phi) -            \Psi^{(2)}((1-\mu)\phi) ),
\\
\kappa_{\mu,\mu,\phi}        &=&  \phi^2 (\mu   \Psi^{(2)}(\mu\phi) + (1-\mu)  \Psi^{(2)}((1-\mu)\phi) ),
\\
\kappa_{\mu,\phi,\phi}       &=&  \phi   (\mu^2 \Psi^{(2)}(\mu\phi) - (1-\mu)^2\Psi^{(2)}((1-\mu)\phi) ),
\\
\kappa_{\phi,\phi,\phi}      &=&          \mu^3 \Psi^{(2)}(\mu\phi) + (1-\mu)^3\Psi^{(2)}((1-\mu)\phi) - \Psi^{(2)}(\phi),
\\
\kappa_{\mu,\mu,\mu,\mu}     &=&  \phi^4 (      \Psi^{(3)}(\mu\phi) +            \Psi^{(3)}((1-\mu)\phi) ),
\\
\kappa_{\mu,\mu,\mu,\phi}    &=&  \phi^3 (\mu   \Psi^{(3)}(\mu\phi) - (1-\mu)  \Psi^{(3)}((1-\mu)\phi) ),
\\
\kappa_{\mu,\mu,\phi,\phi}   &=&  \phi^2 (\mu^2 \Psi^{(3)}(\mu\phi) + (1-\mu)^2\Psi^{(3)}((1-\mu)\phi) ),
\\
\kappa_{\mu,\phi,\phi,\phi}  &=&  \phi   (\mu^3 \Psi^{(3)}(\mu\phi) - (1-\mu)^3\Psi^{(3)}((1-\mu)\phi) ),
\\
\kappa_{\phi,\phi,\phi,\phi} &=&           \mu^4 \Psi^{(3)}(\mu\phi) + (1-\mu)^4\Psi^{(3)}((1-\mu)\phi) - \Psi^{(3)}(\phi),
\\
\kappa_{\mu,\mu\phi}         &=&  \phi   (       \Psi^{(1)}(\mu\phi) + \Psi^{(1)}((1-\mu)\phi)),
\\
\kappa_{\phi,\mu\phi}        &=&   \mu           \Psi^{(1)}(\mu\phi) - (1-\mu)   \Psi^{(1)}((1-\mu)\phi),
\\
\kappa_{\mu,\mu\mu}          &=& \kappa_{\phi,\mu\mu} = \kappa_{\mu,\phi\phi} = \kappa_{\phi,\phi\phi} = 0.
\end{eqnarray*}
\end{appendices}

\end{document}